\documentclass{amsart}

\usepackage[all]{xy}
\usepackage{tkz-graph}
\usetikzlibrary{arrows}
\usetikzlibrary{positioning}
\usepackage{amssymb,verbatim,epstopdf}

 \newtheorem{thm}{Theorem}[section]

 \newtheorem{prop}[thm]{Proposition}
 \theoremstyle{definition}
 \newtheorem{defn}[thm]{Definition}
 \theoremstyle{remark}
 \newtheorem{rem}[thm]{Remark}
 \newtheorem{ex}{Example}

 \numberwithin{equation}{section}

\newcommand\RR{{\mathbb{R}}}
\newcommand\ZZ{{\mathbb{Z}}}
\newcommand\QQ{{\mathbb{Q}}}
\newcommand\CC{{\mathbb{C}}}

\newcommand\Sym{{\mathsf{Sym}}}
\newcommand\Aut{{\mathsf{Aut}}}

\newcommand\St{{\mathsf{St}}}

\newcommand\Sp{{\mathsf{Sp}}}
\newcommand\Stab{{\mathsf{Stab}}}
\newcommand\boundary{{X^\omega}}
\newcommand\mmod{{\mathsf{mod}}}

\newcommand\GG{{\mathcal{G}}}
\newcommand\HH{{\mathcal{H}}}
\newcommand\TT{{\mathcal{T}}}
\newcommand\BB{{\mathcal{B}}}
\newcommand\LL{{\mathcal{L}}}
\newcommand\MM{{\mathcal{M}}}
\newcommand\WW{{\mathcal{W}}}
\newcommand\unitary{{\mathcal{U}}}
\newcommand\yy{{\widetilde{y}}}

\newcommand\fin{{\mathsf{fin}}}

\newcommand{\img}[1]{{\mathsf{IMG}\left(#1\right)}}
\newcommand{\sset}[2]{{\left\{\ #1 \ \mid \ #2 \ \right\}}}

\newcommand\qqand{\qquad\text{and}\qquad}

\begin{document}

\title[From self-similar groups to self-similar sets and spectra]
 {From self-similar groups \\ to self-similar sets and spectra}

\author[Rostislav Grigorchuk]{Rostislav Grigorchuk}
\email{grigorch@math.tamu.edu}
\thanks{The authors acknowledge partial support by the US National Science
Foundation under Grants No. DMS-1105520 and DMS-1207699}

\author[Volodymyr Nekrashevych]{Volodymyr Nekrashevych}
\email{nekrash@math.tamu.edu}

\author[Zoran \v{S}uni\'c]{Zoran \v{S}uni\'c}
\email{sunic@math.tamu.edu}

\address{\\
Department of Mathematics \\
Texas A\&M University \\
College Station, TX 77843-3368, USA}

\subjclass[2010]{Primary 37A30; Secondary 28A80; 05C50, 20E08, 43A07}

\keywords{self-similar groups, Schreier graphs, spectra, self-similar sets, fractals, Julia sets,
amenable action, Laplacian}

\dedicatory{To the memory of Gilbert Baumslag, a great colleague and a great friend}

\begin{abstract}
The survey presents developments in the theory of self-similar groups leading to applications to
the study of fractal sets and graphs and their associated spectra.
\end{abstract}

\maketitle

\section{Introduction}

The purpose of this survey is to present some recent developments in the theory of self-similar
groups and its applications to the study of fractal sets. For brevity, we will concentrate only on
the following two aspects (for other aspects see~\cite{bartholdi-g-n:fractal}):

\begin{enumerate}
\item[(i)] Construction of new fractals by using algebraic tools and interpretation of well known
    fractals (the first Julia set, Sierpi\'nski gasket, Basilica fractal, and other Julia sets of
    post-critically finite rational maps on the Riemann sphere) in terms of self-similar groups
    and their associated objects -- Schreier graphs.

\item[(ii)] Study of the spectra of the Laplacian on Schreier graphs of self-similar groups and
    on the associated fractals by appropriate limiting processes.
\end{enumerate}

The presentation will be focused on a few representative examples for which the ``entire program''
(going from a self-similar group to its associated self-similar objects and calculation/description
of their spectra) is successfully implemented, such as the first Grigorchuk group\footnote{the
second and the third author insist on the use of this terminology} $\GG$, the lamplighter group
$\LL_2$, the 3-peg Hanoi Towers group $\HH$, and the tangled odometers group $\TT$, but also some
examples with only partial implementation, such as the Basilica group $\BB$ and the iterated
monodromy group $\img{z^2+i}$.


\section{Self-similar groups and their Schreier graphs}

\subsection{Schreier graphs}

Let $G$ be a finitely generated group, generated by a finite symmetric set $S$ ($S$ being symmetric
means $S=S^{-1}$) acting on a set $Y$ (all actions in this survey will be left actions). The
\emph{Schreier graph} of the action of $G$ on $Y$ with respect to $S$ is the oriented graph
$\Gamma(G,S,Y)$ defined as follows. The vertex set of the Schreier graph is $Y$ and the edge set is
$S \times Y$. For $s \in S$ and $y \in Y$, the edge $(s,y)$ connects $y$ to $sy$. When the graph is
drawn, the edge $(s,y)$ is usually labeled just by $s$, since its orientation from $y$ to $sy$
uniquely indicates the correct ``full label'' $(s,y)$. In other words, one usually draws
 $y~\bullet \xrightarrow[\phantom{weeert}]{s}     \bullet~sy$ instead of
 $y~\bullet \xrightarrow[\phantom{weeert}]{(s,y)} \bullet~sy$.
 
The Schreier graph $\Gamma(G,S,Y)$ is connected if and only if the action is transitive (some
authors define Schreier graphs only in the transitive/connected case).

\begin{ex}
Let $Y=\{1,2,3,4\}$ and $D_4$ be the subgroup of the symmetric group on $Y$ (with its usual left
action) generated by $S = \langle \sigma,\bar{\sigma},\tau\rangle$, where $\sigma$ is the 4-cycle
$\sigma=(1234)$, $\bar{\sigma}$ is its inverse $\bar{\sigma}= \sigma^{-1} = (1432)$, and $\tau$ is
the transposition $\tau=(24)$ (note that one can interpret $D_4$ as the dihedral group of
isometries of a square with vertices 1,2,3,4; $\sigma$ is the rotation by $\pi/2$ and $\tau$ the
mirror symmetry with respect to the line 13). The Schreier graph $\Gamma(D_4,S,Y)$ is drawn on the
left in Figure~\ref{f:d4}.
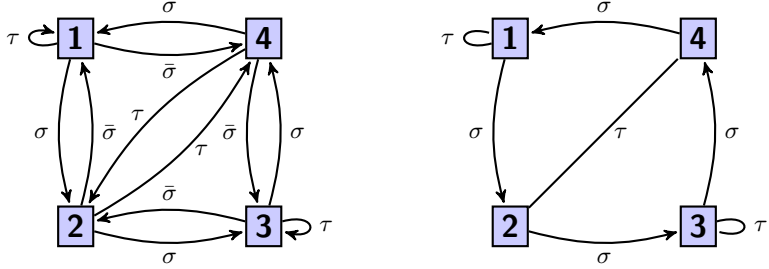
\begin{figure}[!ht]
\begin{tabular}{ccc}
\begin{tikzpicture}[->,>=stealth',shorten >=1pt,auto,node distance=2.5cm,
  thick,main node/.style={rectangle,fill=blue!20,draw,font=\sffamily\Large\bfseries,label distance=-1mm}]

  \node[main node] (1)              {1};
  \node[main node] (2) [below of=1] {2};
  \node[main node] (4) [right of=1] {4};
  \node[main node] (3) [right of=2] {3};

  \path[every node/.style={font=\sffamily\small}]
    (1) edge [bend right=15] node[left]  {$\sigma$}       (2)
        edge [bend right=15] node[below] {$\bar{\sigma}$} (4)
        edge [loop left]     node        {$\tau$}         (1)
    (2) edge [bend right=15] node[below] {$\sigma$}       (3)
        edge [bend right=15] node[right] {$\bar{\sigma}$} (1)
        edge [bend right=15] node[right] {$\tau$}         (4)
    (3) edge [bend right=15] node[right] {$\sigma$}       (4)
        edge [bend right=15] node[above] {$\bar{\sigma}$} (2)
        edge [loop right]    node        {$\tau$}         (3)
    (4) edge [bend right=15] node[above] {$\sigma$}       (1)
        edge [bend right=15] node[left]  {$\bar{\sigma}$} (3)
        edge [bend right=15] node[left]  {$\tau$}         (2);
\end{tikzpicture}
& \hspace{0.5cm} &
\begin{tikzpicture}[->,>=stealth',shorten >=1pt,auto,node distance=2.5cm,
  thick,main node/.style={rectangle,fill=blue!20,draw,font=\sffamily\Large\bfseries}, every loop/.style={-}]

  \node[main node] (1)              {1};
  \node[main node] (2) [below of=1] {2};
  \node[main node] (4) [right of=1] {4};
  \node[main node] (3) [right of=2] {3};

  \path[every node/.style={font=\sffamily\small}]
    (1) edge[bend right=15] node[left]  {$\sigma$} (2)
        edge[loop left]     node        {$\tau$}   (1)
    (2) edge[bend right=15] node[below] {$\sigma$} (3)
        edge[style={-}]     node[right] {$\tau$}   (4)
    (3) edge[bend right=15] node[right] {$\sigma$} (4)
        edge[loop right]    node        {$\tau$}   (3)
    (4) edge[bend right=15] node[above] {$\sigma$} (1);
\end{tikzpicture}
\end{tabular}
\caption{The Schreier graph $\Gamma(D_4,S,Y)$, and its simplified drawing}
\label{f:d4}
\end{figure}
\end{ex}

The edge $(s,y)$ connects $y$ to $sy$ and the edge $(s^{-1},sy)$ goes in the opposite direction and
connects $sy$ to $y$. In order to avoid clutter in the drawings, for each pair of mutually inverse
generators $s,s^{-1}\in S$ that are not involutions, one usually chooses one of them, say $s$, and
only draws the oriented edges labeled by $s$, while all edges labeled by $s^{-1}$ are suppressed.
Further, for an involution $s \in S$ and $y \in Y$, only one unoriented edge is drawn between $y$
and $sy$ (see the graph on the right in Figure~\ref{f:d4} and note that $\sigma$ is not an
involution, while $\tau$ is).

\subsection{Random walk operators on Schreier graphs}

The Schreier graph $\Gamma=\Gamma(G,S,Y)$ is regular with every vertex having both the out-degree
and the in-degree equal to $|S|$. The \emph{random walk} operator on $\Gamma$ (also known as the
\emph{Markov operator}) is the operator
\begin{gather*}
 M : \ell^2(\Gamma) \to \ell^2(\Gamma) \\
 (Mf)(y) = \frac{1}{|S|} \sum_{s \in S} f(sy).
\end{gather*}
where $\ell^2(\Gamma)=\ell^2(Y)$ is the Hilbert space of square summable functions on $Y$
\[
 \ell^2(\Gamma) = \ell^2(Y)  = \sset{f: Y \to \RR}{ \sum_{y \in Y} |f(y)|^2 < \infty}.
\]
Thus, given a function $f: Y \to \RR$ on the vertex set $Y$, the operator $M$ produces an updated
function $Mf:Y \to \RR$ by replacing the value at each vertex $y$ by the average of the $f$-values
at the neighbors of $y$ in the Schreier graph.

For $x \in \RR$, let $M(x)$ be the operator $M(x) = M-xI$. The \emph{spectrum} $\Sp(M)$ of $M$ is
the set of values of $x$ for which the operator $M(x)$ from the pencil of operators $\{ M(x) \mid x
\in \RR \}$ is not invertible. Note that the operator $M$ is bounded (in fact $||M||\leq 1$) and,
since $S$ is symmetric, it is self-adjoint. Therefore its spectrum is a closed subset of the
interval $[-1,1]$. When $Y$ is finite, the spectrum $\Sp(M)$ is just the set of eigenvalues of the
operator $M$, but in general the spectrum only contains the set of eigenvalues of $M$. Recall that
$\lambda$ is an \emph{eigenvalue} of $M$ if and only if $Mf = \lambda f$, for some nonzero function
$f \in \ell^2(\Gamma)$; such a nonzero function is called an \emph{eigenfunction} of $M$.

Let $G=\langle S \rangle$ act on two sets $Y$ and $\widetilde{Y}$ and $\delta: \widetilde{Y} \to Y$
be a surjective $G$-equivariant map, that is, a surjective function $\delta$ such that
$g\delta(\yy) = \delta(g\yy)$, for $g \in G$ and $\yy \in \widetilde{Y}$ (equivalently,
$s\delta(\yy) = \delta(s\yy)$, for $s \in S$ and $\yy \in \widetilde{Y}$). On the level of Schreier
graphs $\delta$ induces a surjective graph homomorphism from
$\Gamma_{\widetilde{Y}}=\Gamma(G,S,\widetilde{Y})$ to $\Gamma_Y=\Gamma(G,S,Y)$ preserving edge
labels and sending the edge $\yy \stackrel{s}{\to} s\yy$ to the edge $\delta(\yy) \stackrel{s}{\to}
s\delta(\yy)$. We say that $\Gamma_{\widetilde{Y}}$ is a \emph{covering} of $\Gamma_{Y}$ and
$\delta$ is a \emph{covering map}.

Assume that both $\widetilde{Y}$ and $Y$ are finite. For every function $f \in \ell^2(\Gamma_Y)$,
define the \emph{lift} $\widetilde{f} \in \ell^2(\Gamma_{\widetilde{Y}})$ by $\widetilde{f}(\yy) =
f(\delta(\yy))$, for $\yy \in \widetilde{Y}$. For all $f \in \ell^2(\Gamma_Y)$, we have
\[
 (M_{\widetilde{Y}} \widetilde{f}) (\yy) = (M_Yf)(\delta(\yy)).
\]
If $f$ is an eigenfunction of $M_Y$ with eigenvalue $\lambda$, then $\widetilde{f}$ is an
eigenfunction of $M_{\widetilde{Y}}$ with the same eigenvalue. Therefore, whenever there exists a
surjective $G$-equivariant map $\delta:\widetilde{Y} \to Y$ between two finite sets $\widetilde{Y}$
and $Y$, the spectrum of $M_Y$ is included in the spectrum of $M_{\widetilde{Y}}$, that is,
$\Sp(M_Y) \subseteq \Sp(M_{\widetilde{Y}})$.

Let $\{Y_n\}_{n=0}^\infty$ be a sequence of finite $G$-sets (sets with a $G$-action defined on
them), $\{\delta_n:Y_{n+1} \to Y_n\}_{n=0}^\infty$ a sequence of surjective $G$-equivariant maps,
$Y$ be a $G$-set, and $\{\widetilde{\delta}_n:Y \to Y_n\}_{n=0}^\infty$ a sequence of surjective
$G$-equivariant maps such that $\delta_n\widetilde{\delta}_{n+1}=\widetilde{\delta}_n$, for $n \geq
0$. Denote $\Gamma_n = \Gamma(G,S,Y_n)$, $\Gamma=\Gamma(G,S,Y)$, and the corresponding random walk
operators by $M_n$ and $M$, respectively. The sequences of equivariant maps $\{\delta_n\}$ and
$\{\widetilde{\delta}_n\}$ induce graph coverings between the corresponding Schreier graphs such
that the following diagram commutes
\begin{equation}\label{d:equiv}
\xymatrix{
 \Gamma_0 \ar@{<-}[r]^{\delta_0} & \Gamma_1 \ar@{<-}[r]^{\delta_1}  & \Gamma_2 \ar@{<-}[r]^{\delta_2} & \dots \\
 \Gamma   \ar@{->}[u]^{\widetilde{\delta}_0} \ar@{->}[ur]^{\widetilde{\delta}_1} \ar@{->}[urr]^{\widetilde{\delta}_2} \ar@{..>}[urrr]
}
\end{equation}
and we obtain an increasing sequence $\{\Sp(M_n)\}_{n=0}^\infty$ of finite sets, each consisting of
the eigenvalues of $M_n$. We are interested in situations in which this sequence is sufficient to
determine the spectrum of $M$ in the sense that
\[
 \overline{\bigcup_{n=0}^\infty \Sp(M_n)} = \Sp(M).
\]

\begin{ex}\label{ex:dinf}
This example is relatively straightforward, but it illustrates the setup we introduced above.
Consider the infinite dihedral group $D_\infty = \langle a,b\rangle$, generated by two involutions
$a$ and $b$. We may think of it as the group of isometries of the set of integer points on the real
line, with the action of $a$ and $b$ given by $a(n)=1-n$ and $b(n)=-n$. Let $Y= \ZZ$ and $\Gamma$
be the Schreier graph $\Gamma=\Gamma(D_\infty,S,Y)$, drawn in the bottom row in
Figure~\ref{f:dinf}. For $n \geq 0$, let $Y_n=\{0,\pm1,\dots,\pm 2^{n-1}-1, 2^{n-1}\}$. Note that
$Y_n$ is a set of unique representatives of the residue classes modulo $2^n$, for $n \geq 0.$  Thus
we may think of $Y_n$ as $\ZZ/2^n\ZZ$. The action of $D_\infty$ on $\ZZ$ induces a well defined
action on the set of residue classes $\ZZ/2^n\ZZ$, for $n \geq 0$, and we denote
$\Gamma_n=\Gamma(D_\infty,S,Y_n).$ The sequence of Schreier graphs $\{\Gamma_n\}$ is indicated in
the top row in Figure~\ref{f:dinf}. For $n \geq 0$, the maps $\delta_n: Y_{n+1} \to Y_n$ and
$\widetilde{\delta}_n: Y \to Y_n$, given by $\delta_n(y) = \mmod(y,2^n)$, for $y \in
\ZZ/2^{n+1}\ZZ$, and $\widetilde{\delta}_n(y) = \mmod(y,2^n)$, for $y \in \ZZ$, where
$\mmod(y,2^n)$ is the remainder obtained when $y$ is divided by $2^n$, are $D_\infty$-equivariant.
\begin{figure}[!ht]
\begin{tikzpicture}[-,>=stealth',shorten >=1pt,auto,node distance=1.5cm,thick,
  main node/.style={rectangle,fill=blue!20,draw,font=\sffamily\Large\bfseries},
  dots node/.style={},
  every loop/.style={-}]

  \node[main node] (1)              {0};
  \node[main node] (2) [right of=1] {0};
  \node[main node] (3) [right of=2] {1};
  \node[main node] (4) [right of=3] {0};
  \node[main node] (5) [right of=4] {1};
  \node[main node] (6) [right of=5] {-1};
  \node[main node] (7) [right of=6] {2};
  \node[dots node] (8) [right of=7] {...};
  \node[main node] (9) [below=0.7cm of 1] {0};
  \node[main node] (10)[below=0.7cm of 2] {1};
  \node[main node] (11)[below=0.7cm of 3] {-1};
  \node[main node] (12)[below=0.7cm of 4] {2};
  \node[main node] (13)[below=0.7cm of 5] {-2};
  \node[main node] (14)[below=0.7cm of 6] {3};
  \node[main node] (15)[below=0.7cm of 7] {-3};
  \node[dots node] (16)[below=0.95cm of 8] {...};

  \path[every node/.style={font=\sffamily\small}]
    (1) edge[loop right] node[above] {$a$}  (1)
        edge[loop left]  node[above] {$b$}  (1)
    (2) edge[loop left]  node[above] {$b$}  (3)
        edge             node[above] {$a$}  (3)
    (3) edge[loop right] node[above] {$b$}  (3)
    (4) edge[loop left]  node[above] {$b$}  (4)
        edge             node[above] {$a$}  (5)
    (5) edge             node[above] {$b$}  (6)
    (6) edge             node[above] {$a$}  (7)
    (7) edge[loop right] node[above] {$b$}  (7)
    (9) edge[loop left]  node[above] {$b$}  (9)
    (9) edge             node[above] {$a$}  (10)
    (10)edge             node[above] {$b$}  (11)
    (11)edge             node[above] {$a$}  (12)
    (12)edge             node[above] {$b$}  (13)
    (13)edge             node[above] {$a$}  (14)
    (14)edge             node[above] {$b$}  (15)
    (15)edge             node[above]  {$a$} (16)
    ;
\end{tikzpicture}
\caption{The Schreier graphs $\Gamma_0,\Gamma_1,\Gamma_2,\dots$ (top row) and $\Gamma$ (bottom row) for the infinite dihedral group $D_\infty=\langle a,b\rangle$}
\label{f:dinf}
\end{figure}
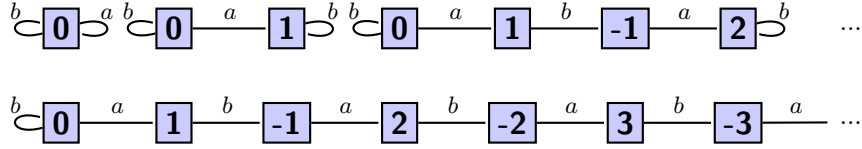

For $n \geq 0$, $\Sp(M_n)$ consists of $2^n$ distinct eigenvalues of multiplicity 1
\[
 \Sp(M_n) = \{1\} \cup \frac{1}{2} \bigcup_{i=0}^{n-1} f^{-i}(0) =
 \{1,0\} \cup \frac{1}{2}\scriptstyle{\sset{\underbrace{\pm\sqrt{2 \pm \sqrt{2 \pm \dots \pm \sqrt{2}}}}_{i \textup{ plus-minus signs}}}{i=1,\dots,n-1}},
\]
where $f(x)=x^2-2$. On the other hand, the spectrum of the doubly infinite path $\Gamma$ is $[-1,1]$ and we have
\[
 \Sp(M) = [-1,1] = \frac{1}{2}\overline{\bigcup_{n=0}^{\infty} f^{-n}(0)} = \overline{\bigcup_{n=0}^\infty \Sp(M_n)}.
\]
\end{ex}

\subsection{Adjacency operator on Schreier graphs and Schreier spectrum}

For the Schreier graph $\Gamma=\Gamma(G,S,Y)$ of the action of $G=\langle S \rangle$ on $Y$, the
\emph{adjacency operator} on $\Gamma$ is the operator $A: \ell^2(\Gamma) \to \ell^2(\Gamma)$
defined by
\[
 (Af)(y) = \sum_{s \in S} f(sy).
\]
The random walk operator $M=\frac{1}{|S|}A$ is the normalized version of the adjacency operator $A$
and their spectra are just multiples of each other. Denote the spectrum of $A$ by $\Sp(\Gamma)$ and
call it the \emph{Schreier spectrum} of $\Gamma$. This is the so called \emph{adjacency spectrum},
but we want to emphasize the scope of all our considerations, namely, adjacency spectra of Schreier
graphs of finitely generated groups. For the purposes of our calculations, the Schreier spectra
turn out to be the most convenient, but it is easy to switch to their Markovian or Laplacian
versions when needed (the \emph{Laplacian} operator is the operator $L=I-M$, where $I$ is the
identity operator).

\subsection{Rooted regular trees and self-similar groups}

We introduce the class of self-similar groups acting on regular rooted trees, providing a framework
for examples like Example~\ref{ex:dinf}, and a source of other examples.

Let $X$ be a finite set, usually called the \emph{alphabet}, of size $k$. The set of all finite
words over $X$ is denoted by $X^*$. The set $X^*$ can be naturally equipped with the structure of a
\emph{rooted $k$-regular tree} as follows. The vertices of the tree are the words in $X^*$, the
\emph{root} is the empty word $\epsilon$, the \emph{level} $n$ is the set $X^n$ of words of length
$n$ over $X$, and the children of each vertex $u \in X^*$ are the $k$ vertices of the form $ux$,
for $x \in X$. We use $X^*$ to denote the set of finite words over $X$, the set of vertices of the
rooted tree we just described, as well as the tree itself.

The group $\Aut(X^*)$ of all automorphisms of the rooted $k$-regular tree $X^*$ preserves the root
and all levels of the tree. Every automorphism $g \in \Aut(X^*)$ induces a permutation $\alpha_g$
of $X$, defined by $\alpha_g(x) = g(x)$, called the \emph{root permutation} of $g$. It represents
the action of $g$ at the first letter in each word. For every automorphism $g \in \Aut(X^*)$ and
every vertex $u \in X^*$, there exists a unique tree automorphism of $X^*$, denoted by $g_u$, such
that, for all words $w \in X^*$,
\[
 g(uw) = g(u)g_u(w).
\]
The automorphism $g_u$ is called the \emph{section} of $g$ at $u.$  It represents the action of
$g$ on the tails of words that start with $u$. Every automorphism $g$ is uniquely determined by its
root permutation $\alpha_g$ and the $k$ sections at the first level $g_x$, for $x \in X$. Indeed,
for every $x \in X$ and $w \in X^*$ we have
\begin{equation}\label{e:gxw}
 g(xw) = \alpha_g(x) g_x(w).
\end{equation}

When $X=\{0,1,\dots,k-1\}$, a succinct representation, called \emph{wreath recursion}, of the
automorphism $g \in \Aut(X^*)$, describing its root permutation and its first level sections is
given by
\begin{equation}\label{e:wreath}
 g = \alpha_g(g_0,g_1,\dots,g_{k-1}).
\end{equation}
In addition of being short and clear, it has many other advantages, not the least of which is that
it emphasizes the fact that $\Aut(X^*)$ is isomorphic to the semidirect product $\Sym(X) \ltimes
(\Aut(X^*))^X$, that is, to the permutational wreath product $\Sym(X) \wr_X \Aut(X^*)$, where
$\Sym(X)$ is the group of all permutations of $X$.

A set $S \subseteq \Aut(X^*)$ of tree automorphisms is \emph{self-similar} if it is closed under
taking sections, that is, every section of every element of $S$ is itself in the set $S$. Thus, for
every word $u$, the action of every automorphism $s \in S$ on the tails of words that start with
$u$ looks exactly like the action of some element of $S$. Note that for a set $S$ to be
self-similar it is sufficient that it contains the first level sections of all of its elements.
Indeed, this is because $g_{uv} = (g_u)_v$, for all words $u,v \in X^*$. A group $G \leq \Aut(X^*)$
of tree automorphisms is \emph{self-similar} if it is self-similar as a set. Every group generated
by a self-similar set is itself self-similar. This is because ``sections of the product are
products of sections'' and ``sections of the inverse are inverses of sections''. To be precise, for
all tree automorphisms $g$ and $h$ and all words $u \in X^*$,
\[
 (gh)_u = g_{h(u)}h_u \qqand \left(g^{-1}\right)_u = \left(g_{g^{-1}(u)}\right)^{-1}
\]
The observation that groups generated by self-similar sets are themselves self-similar enables one
to easily construct many examples of finitely generated self-similar groups, as demonstrated in the
next subsection.

\begin{rem}
It should be clarified that when we speak of a subset $S$ or a subgroup $G$ of $\Aut(X^*)$ as a
self-similar set, we do not use this terminology in the, by now widely accepted and used, sense of
Hutchinson~\cite{hutchinson:self-similarity}. It would be more precise to say, and it is often
said, that the action is self-similar, that is, the action is adapted to the self-similar nature of
the rooted tree and its boundary, the Cantor set. Self-similar sets in the sense of Hutchinson do
play a role here, as such sets appear as limit spaces of contracting self-similar groups (see
Section~\ref{s:img}) and our considerations lead to results on Laplacians on such self-similar sets
(see Section~\ref{s:laplacian}).
\end{rem}

\subsection{Automaton groups}

An \emph{automaton}, in our context, is any finite self-similar set $S$ of tree automorphisms. The
group $G(S) = \langle S \rangle$, called the \emph{automaton group} over $S$ (or of $S$),  is a
finitely generated self-similar group. A simple way to define an automaton is by defining the
action of each of its elements recursively as in~\eqref{e:gxw}.

\begin{ex}
Consider the binary rooted tree based on the alphabet $X=\{0,1\}^*.$  Define a finite self-similar
set $S=\{a,b\}$ of tree automorphisms recursively by
\begin{alignat*}{2}
 a(0u) &= 1a(u) \qquad b(0u) &&= 0b(u), \\
 a(1u) &= 0b(u) \qquad b(1u) &&= 1a(u),
\end{alignat*}
for every word $u \in X^*$, and $a(\epsilon)=b(\epsilon) = \epsilon$. Evidently, the root
permutations and the sections of $a$ and $b$ are given in the following table.
\[
 \begin{array}{c||c|cc}
 s & \alpha_s & s_0 & s_1 \\
 \hline
 a & (01) & a & b \\
 b & ()   & b & a
 \end{array}
\]
where $()$ and $(01)$ denote, respectively, the trivial and the nontrivial permutation of
$X=\{0,1\}$. Calculating the action of any element of $S$ on any word in $X^*$ by using the
recursive definition is straightforward. For instance,
\[
 a(10101) = 0b(0101) = 00b(101) = 001a(01) = 0011a(1) = 00110.
\]
\end{ex}

One may think of the elements of an automaton $S$ as the \emph{states} of a certain type of
transducer, a so-called \emph{Mealy machine}. The recursive definition~\ref{e:gxw} of the action
of $s \in S$ is interpreted as follows. To calculate the action of the state $s$ on some input
word $xu$ starting with $x$, the machine first rewrites $x$ into $\alpha_s(x)$, changes its state
to $s_x$, and lets the new state handle the rest of the input $u$ in the same manner.  It reads the
first letter of $u$, rewrites it appropriately, then moves to an appropriate state, which then
handles the rest of the input, and so on, until the entire input word is read. It is common to
represent the automaton $S$ by an oriented labeled graph as follows. The vertex
set is the set of states $S$, and each pair of a state $s \in S$ and a letter $x \in X$ determines
a directed edge from $s$ to $s_x$ labeled by $x|\alpha_s(x)$ (equivalently, by $s|s(x)$).

\begin{ex}\label{ex:automata}
Four examples of finite self-similar sets of tree automorphisms are given in Figure~\ref{f:4ex}.
The self-similar groups defined by these sets are the lamplighter group $\LL_2 = \ZZ \ltimes \left(
\oplus_{\ZZ} \ZZ/2\ZZ \right)$ (top left), the dihedral group $D_\infty$ (top right), the binary
odometer group $\ZZ$ (bottom left), and the tangled odometers group $\TT$ (bottom right).
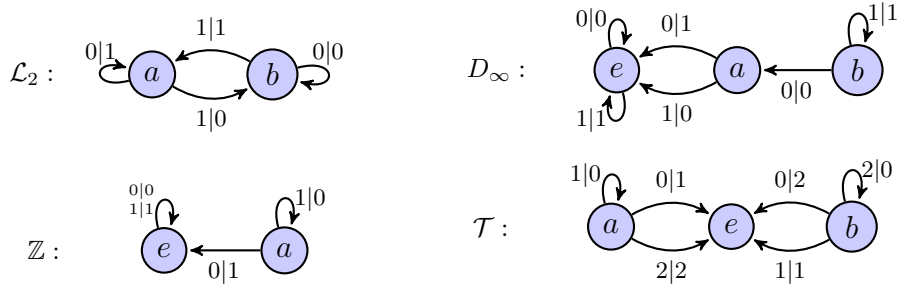
\begin{figure}[!ht]
\begin{tabular}{ccc}
\begin{tikzpicture}[->,>=stealth',shorten >=1pt,auto,node distance=1.6cm,thick,
  main node/.style={circle,fill=blue!20,draw,font=\sffamily\Large\bfseries},
  dots node/.style={}]

  \node[main node] (a)              {$a$};
  \node[main node] (b) [right of=a] {$b$};
  \node[dots node] (0) [left  of=a] {$\LL_2:$};

  \path[every node/.style={font=\sffamily\small}]
    (a) edge[loop left]  node[above] {$0|1$} (a)
        edge[bend right] node[below] {$1|0$} (b)
    (b) edge[loop right] node[above] {$0|0$} (b)
        edge[bend right] node[above] {$1|1$} (a)
    ;
\end{tikzpicture}
& \hspace{7mm}
&
\begin{tikzpicture}[->,>=stealth',shorten >=1pt,auto,node distance=1.6cm,thick,
  main node/.style={circle,fill=blue!20,draw,font=\sffamily\Large\bfseries},
  dots node/.style={}]

  \node[main node] (e)              {$e$};
  \node[main node] (a) [right of=e] {$a$};
  \node[main node] (b) [right of=a] {$b$};
  \node[dots node] (0) [left  of=e] {$D_\infty:$};

  \path[every node/.style={font=\sffamily\small}]
    (e) edge[loop above] node[left]  {$0|0$} (e)
        edge[loop below] node[left]  {$1|1$} (e)
    (a) edge[bend right] node[above] {$0|1$} (e)
        edge[bend left]  node[below] {$1|0$} (e)
    (b) edge             node        {$0|0$} (a)
        edge[loop above] node[right] {$1|1$} (b)
    ;
\end{tikzpicture}
\\
\begin{tikzpicture}[->,>=stealth',shorten >=1pt,auto,node distance=1.6cm,thick,
  main node/.style={circle,fill=blue!20,draw,font=\sffamily\Large\bfseries},
  dots node/.style={}]

  \node[main node] (e)              {$e$};
  \node[main node] (a) [right of=e] {$a$};
  \node[dots node] (0) [left  of=e] {$\ZZ:$};

  \path[every node/.style={font=\sffamily\small}]
    (e) edge[loop above] node[left]  {$\substack{0|0\\1|1}$} (e)
    (a) edge             node        {$0|1$}                 (e)
        edge[loop above] node[right] {$1|0$}                 (a)
    ;
\end{tikzpicture}
& \qquad  &
\begin{tikzpicture}[->,>=stealth',shorten >=1pt,auto,node distance=1.6cm,thick,
  main node/.style={circle,fill=blue!20,draw,font=\sffamily\Large\bfseries},
  dots node/.style={}]

  \node[main node] (a)              {$a$};
  \node[main node] (e) [right of=a] {$e$};
  \node[main node] (b) [right of=e] {$b$};
  \node[dots node] (0) [left  of=a] {$\TT:$};

  \path[every node/.style={font=\sffamily\small}]
    (a) edge[loop above] node[left]  {$1|0$} (a)
        edge[bend left ] node[above] {$0|1$} (e)
        edge[bend right] node[below] {$2|2$} (e)
    (b) edge[bend right] node[above] {$0|2$} (e)
        edge[bend left]  node[below] {$1|1$} (e)
        edge[loop above] node[right] {$2|0$} (b)
    ;
\end{tikzpicture}
\end{tabular}
\caption{Automata defining $\LL_2$, $D_\infty$, $\ZZ$, and $\TT$}
\label{f:4ex}
\end{figure}
In the last three automata the state $e$ represents the trivial automorphism of the tree, which
does not change any input word. Thus, we use $\epsilon$ for the empty word, that is, the root of
$X^*$, $()$ for the trivial permutation of $X$, and $e$ for the trivial automorphism of the tree
$X^*$. To avoid clutter, in the automaton for $\ZZ$ we used the convention that the same edge may
be used with several labels, while in the automaton for $\TT$ the convention that the loops
associated to the trivial state $e$ are not drawn. Note that the first three automata are defined
over the binary alphabet $X=\{0,1\}$ while the last one is defined over the ternary alphabet
$X=\{0,1,2\},$ hence that group acts on the ternary rooted tree.

One can easily switch back and forth between the various representations of the given automata. For
instance, the recursive definition of the action of the dihedral group $D_\infty = \langle
a,b\rangle$ on the binary rooted tree is given by
\begin{align*}
 a(0u) &= 1u, \qquad & b(0u) = 0a(u), \\
 a(1u) &= 0u, \qquad & b(1u) = 1b(u),
\end{align*}
Tabular representation of the self-similar set defining $\TT$ and the wreath recursion describing
the same set are given on the left and on the right, respectively in
\begin{equation}\label{e:tangled}
\begin{tabular}{ccc}
 $\begin{array}{c||c|ccc}
 s & \alpha_s & s_0 & s_1 & s_2\\
 \hline
 a & (01) & e & a & e  \\
 b & (02) & e & e & b
 \end{array}$
& \hspace{2cm} &
 $\begin{array}{l}
  \\
  a = (01)(e,a,e) \\
  b = (02)(e,e,b)
 \end{array}$
\end{tabular}
\end{equation}
\end{ex}

It is clear that defining a finitely generated self-similar group is an easy task, in particular
for automaton groups (note that not all finitely generated self-similar groups are automaton
groups). One can methodically construct, one by one, all automaton groups by constructing all
automata with a given number of states over an alphabet of a given size. However, it is not an easy
task to recognize the group that is generated by a given automaton. A full classification of all
automaton groups defined by automata with given number of states $m$ and size of the alphabet $k$
has been achieved only for $m=k=2$~\cite{grigorchuk-n-s:automata}, while for the next smallest case
$m=3$ and $k=2$ only a partial classification was obtained \cite{bondarenko-al:classification32}.

\subsection{The boundary action and the convergence $\Gamma_n \to \Gamma$}

Let $G=\langle S \rangle$, with $S$ symmetric and finite, be a finitely generated subgroup of
$\Aut(X^*)$ and, for $n \geq 0$, let $\Gamma_n = \Gamma(G,S,X^n)$ be the corresponding Schreier
graph of the action on level $n$. The map $\delta_n:X^{n+1} \to X^n$ given by deleting the last
letter in each word is $G$-equivariant and induces a sequence of coverings of degree $|X|$
\[
\xymatrix{
 \Gamma_0 \ar@{<-}[r]^{\delta_0} & \Gamma_1 \ar@{<-}[r]^{\delta_1}  & \Gamma_2 \ar@{<-}[r]^{\delta_2} & \dots
}
\]
Under the covering $\delta_n$ each of the $|X|$ edges $ux~\bullet \xrightarrow[\phantom{weert}]{s}
\bullet~s(u)s_u(x)$ in $\Gamma_{n+1}$, for $x \in X$, is mapped to the edge $u~\bullet
\xrightarrow[\phantom{weert}]{s} \bullet~s(u)$ in $\Gamma_n$.

\begin{ex}
The \emph{first Grigorchuk group} $\GG$ is the self-similar group $\GG=\langle a,b,c,d \rangle$
generated by four involutions $a$, $b$, $c$, and $d$ acting on the binary tree and given by the
wreath recursion
\[
 a = (01)(e,e), \qquad b=()(a,c), \qquad c=()(a,d), \qquad d=()(e,b).
\]
The Schreier graphs of its action on levels 0,1,2, and 3, are given in Figure~\ref{f:gg}.
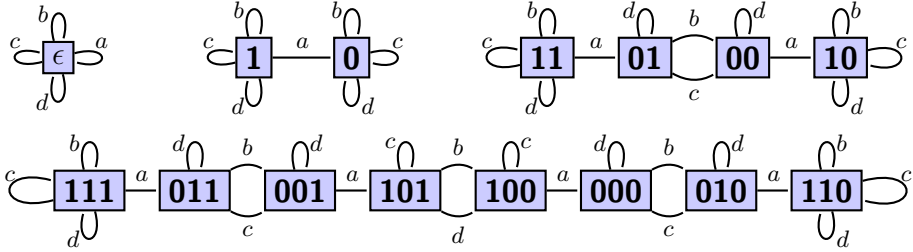
\begin{figure}[!ht]
\begin{tikzpicture}[-,>=stealth',shorten >=1pt,auto,node distance=1.3cm,thick,
  main node/.style={rectangle,fill=blue!20,draw,font=\sffamily\Large\bfseries},
  dots node/.style={},
  every loop/.style={-}]

  \node[main node] (e)                {$\epsilon$};
  \node[dots node] (x)  [right of=e]  {};
  \node[main node] (1)  [right of=x]  {1};
  \node[main node] (0)  [right of=1]  {0};
  \node[dots node] (y)  [right of=0]  {};
  \node[main node] (11) [right of=y]  {11};
  \node[main node] (01) [right of=11] {01};
  \node[main node] (00) [right of=01] {00};
  \node[main node] (10) [right of=00] {10};

  \path[every node/.style={font=\sffamily\small}]
    (e) edge[loop right] node[above] {$a$}  (e)
        edge[loop above] node[left]  {$b$}  (e)
        edge[loop left ] node[above] {$c$}  (e)
        edge[loop below] node[left]  {$d$}  (e)
    (1) edge[loop above] node[left]  {$b$}  (1)
        edge[loop left ] node[above] {$c$}  (1)
        edge[loop below] node[left]  {$d$}  (1)
        edge             node[above] {$a$}  (0)
    (0) edge[loop above] node[left]  {$b$}  (0)
        edge[loop right] node[above] {$c$}  (0)
        edge[loop below] node[right] {$d$}  (0)
    (11)edge[loop above] node[left]  {$b$}  (11)
        edge[loop left ] node[above] {$c$}  (11)
        edge[loop below] node[left]  {$d$}  (11)
        edge             node[above] {$a$}  (01)
    (01)edge[loop above] node[left]  {$d$}  (01)
        edge[bend left]  node[above] {$b$}  (00)
        edge[bend right] node[below] {$c$}  (00)
    (00)edge[loop above] node[right] {$d$}  (00)
        edge             node[above] {$a$}  (10)
    (10)edge[loop above] node[right] {$b$}  (10)
        edge[loop right] node[above] {$c$}  (10)
        edge[loop below] node[right] {$d$}  (10)
    ;
\end{tikzpicture}
\begin{tikzpicture}[-,>=stealth',shorten >=1pt,auto,node distance=1.4cm,thick,
  main node/.style={rectangle,fill=blue!20,draw,font=\sffamily\Large\bfseries},
  dots node/.style={},
  every loop/.style={-}]

  \node[main node] (9)               {111};
  \node[main node] (10)[right of=9]  {011};
  \node[main node] (11)[right of=10] {001};
  \node[main node] (12)[right of=11] {101};
  \node[main node] (13)[right of=12] {100};
  \node[main node] (14)[right of=13] {000};
  \node[main node] (15)[right of=14] {010};
  \node[main node] (16)[right of=15] {110};

  \path[every node/.style={font=\sffamily\small}]
    (9) edge[loop above] node[left ] {$b$}  (9)
        edge[loop left ] node[above] {$c$}  (9)
        edge[loop below] node[left ] {$d$}  (9)
        edge             node[above] {$a$} (10)
    (10)edge[loop above] node[left]  {$d$} (10)
        edge[bend left]  node[above] {$b$} (11)
        edge[bend right] node[below] {$c$} (11)
    (11)edge[loop above] node[right] {$d$} (11)
        edge             node[above] {$a$} (12)
    (12)edge[loop above] node[left]  {$c$} (12)
        edge[bend left]  node[above] {$b$} (13)
        edge[bend right] node[below] {$d$} (13)
    (13)edge[loop above] node[right] {$c$} (13)
        edge             node[above] {$a$} (14)
    (14)edge[loop above] node[left]  {$d$} (14)
        edge[bend left]  node[above] {$b$} (15)
        edge[bend right] node[below] {$c$} (15)
    (15)edge[loop above] node[right] {$d$} (15)
        edge             node[above] {$a$} (16)
    (16)edge[loop above] node[right] {$b$}  (16)
        edge[loop right] node[above] {$c$}  (16)
        edge[loop below] node[right] {$d$}  (16)
    ;
\end{tikzpicture}
\caption{The Schreier graphs $\Gamma_0$, $\Gamma_1$, $\Gamma_2$ (top row), and $\Gamma_3$ (bottom row) for the first Grigorchuk group $\GG$}
\label{f:gg}
\end{figure}
This group was constructed by the first author in~\cite{grigorchuk:burnside} as a particularly
simple example of a finitely generated, infinite 2-group. It was the first example of a group of
intermediate growth and the first example of an amenable group that is not elementary amenable
~\cite{grigorchuk:gdegree}(we will get back to this aspect later).
\end{ex}

\begin{ex}\label{ex:basilica}
The Basilica group is the self-similar group $\BB=\langle a,b \rangle$ generated by the binary tree
automorphisms $a$ and $b$ given by the wreath recursion
\[
 a = (01)(e,b), \qquad b=()(e,a).
\]
The Schreier graphs of its action on levels 0,1,2, and 3, are given in Figure~\ref{f:basilica}.
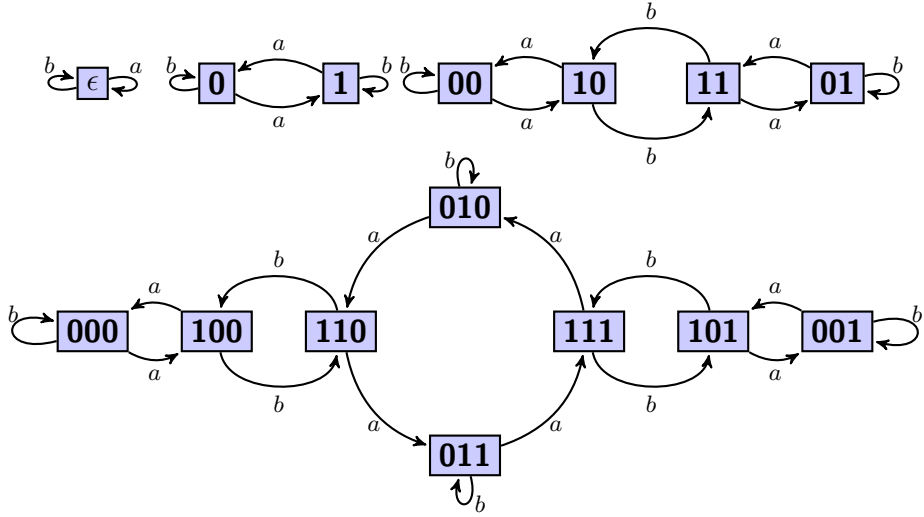
\begin{figure}[!ht]
\begin{tikzpicture}[->,>=stealth',shorten >=1pt,auto,node distance=1.65cm,thick,
  main node/.style={rectangle,fill=blue!20,draw,font=\sffamily\Large\bfseries},
  dots node/.style={}]

  \node[main node] (e)                 {$\epsilon$};
  \node[main node] (0)  [right of=e]   {0};
  \node[main node] (1)  [right of=0]   {1};
  \node[main node] (00) [right of=1]   {00};
  \node[main node] (10) [right of=00]  {10};
  \node[main node] (11) [right of=10]  {11};
  \node[main node] (01) [right of=11]  {01};
  \node[main node] (010)[below of=00]  {010};
  \node[dots node] (c)  [below of=010] {};
  \node[main node] (110)[left  of=c]   {110};
  \node[main node] (100)[left  of=110] {100};
  \node[main node] (000)[left  of=100] {000};
  \node[main node] (111)[right of=c]   {111};
  \node[main node] (101)[right of=111] {101};
  \node[main node] (001)[right of=101] {001};
  \node[main node] (011)[below of=c]   {011};

  \path[every node/.style={font=\sffamily\small}]
    (e) edge[loop right]    node[above] {$a$}  (e)
        edge[loop left]     node[above] {$b$}  (e)
    (0) edge[loop left]     node[above] {$b$}  (0)
        edge[bend right]    node[below] {$a$}  (1)
    (1) edge[loop right]    node[above] {$b$}  (1)
        edge[bend right]    node[above] {$a$}  (0)
   (00) edge[loop left]     node[above] {$b$}  (00)
        edge[bend right]    node[below] {$a$}  (10)
   (10) edge[bend right=80] node[below] {$b$}  (11)
        edge[bend right]    node[above] {$a$}  (00)
   (11) edge[bend right]    node[below] {$a$}  (01)
        edge[bend right=80] node[above] {$b$}  (10)
   (01) edge[loop right]    node[above] {$b$}  (01)
        edge[bend right]    node[above] {$a$}  (11)
  (000) edge[loop left]     node[above] {$b$}  (000)
        edge[bend right]    node[below] {$a$}  (100)
  (100) edge[bend right=80] node[below] {$b$}  (110)
        edge[bend right]    node[above] {$a$}  (000)
  (110) edge[bend right]    node[below] {$a$}  (011)
        edge[bend right=80] node[above] {$b$}  (100)
  (111) edge[bend right=80] node[below] {$b$}  (101)
        edge[bend right]    node[above] {$a$}  (010)
  (101) edge[bend right]    node[below] {$a$}  (001)
        edge[bend right=80] node[above] {$b$}  (111)
  (001) edge[loop right]    node[above] {$b$}  (001)
        edge[bend right]    node[above] {$a$}  (101)
  (010) edge[bend right]    node[above] {$a$}  (110)
        edge[loop above]    node[left]  {$b$}  (010)
  (011) edge[bend right]    node[below] {$a$}  (111)
        edge[loop below]    node[right] {$b$}  (011)
            ;
\end{tikzpicture}
\caption{The Schreier graphs $\Gamma_0$, $\Gamma_1$, $\Gamma_2$ (top row), and $\Gamma_3$ (bottom row) for the Basilica group $\BB$}
\label{f:basilica}
\end{figure}

The group $\BB$ was first considered in~\cite{grigorchuk-z:basilica1}
and~\cite{grigorchuk-z:basilica2} where it was proved that it is a weakly branch, torsion free
group which is not sub-exponentially amenable. It was later proved by Bartholdi and
Vir{\'a}g~\cite{bartholdi-v:basilica}, using speed estimates for random walks, that this group is amenable,
thus providing the first example of an amenable group that is not sub-exponentially amenable.
\end{ex}

\begin{ex}
The Hanoi Towers group is the self-similar group $\HH=\langle a,b,c \rangle$ generated by three
involutions acting on the ternary tree given by the wreath recursion
\[
 a = (01)(e,e,a), \qquad b=(02)(e,b,e), \qquad c = (12)(c,e,e).
\]
The Schreier graphs of its action on levels 0,1, and 2 are given in Figure~\ref{f:hanoi}.
\begin{figure}[!ht]
\begin{tikzpicture}[-,>=stealth',shorten >=1pt,auto,node distance=1.5cm,thick,
  main node/.style={rectangle,fill=blue!20,draw,font=\sffamily\Large\bfseries},
  dots node/.style={}, every loop/.style={-}]

  \node[main node] (e)                      {$\epsilon$};
  \node[main node] (2)  [below right of=e]  {2};
  \node[main node] (0)  [above right of=2]  {0};
  \node[main node] (1)  [below right of=0]  {1};
  \node[main node] (21) [below right of=1]  {21};
  \node[main node] (11) [below left  of=21] {11};
  \node[main node] (01) [below right of=21] {01};
  \node[main node] (20) [above right of=21] {20};
  \node[main node] (00) [above right of=20] {00};
  \node[main node] (10) [below right of=00] {10};
  \node[main node] (12) [below right of=10] {12};
  \node[main node] (22) [below right of=12] {22};
  \node[main node] (02) [below left of=12]  {02};

  \path[every node/.style={font=\sffamily\small}]
    (e) edge[loop right] node[below] {$a$}  (e)
        edge[loop left]  node[below]  {$b$}  (e)
        edge[loop above] node[left]  {$c$}  (e)
    (0) edge[loop above] node[left]  {$c$}  (0)
        edge             node[left]  {$b$}  (2)
        edge             node[right] {$a$}  (1)
    (1) edge[loop right] node[below] {$b$}  (1)
        edge             node[below] {$c$}  (2)
    (2) edge[loop left]  node[below] {$a$}  (2)
   (00) edge[loop above] node[left]  {$c$}  (0)
        edge             node[left]  {$b$}  (20)
        edge             node[right] {$a$}  (10)
   (10) edge             node[right] {$b$}  (12)
        edge             node[below] {$c$}  (20)
   (20) edge             node[left]  {$a$}  (21)
   (21) edge             node[left]  {$c$}  (11)
        edge             node[right] {$b$}  (01)
   (12) edge             node[left]  {$a$}  (02)
        edge             node[right] {$c$}  (22)
   (11) edge[loop left]  node[below] {$b$}  (11)
        edge             node[below] {$a$}  (01)
   (01) edge             node[below] {$c$}  (02)
   (02) edge             node[below] {$b$}  (22)
   (22) edge[loop right] node[below] {$a$}  (22)
            ;
\end{tikzpicture}
\caption{The Schreier graphs $\Gamma_0$, $\Gamma_1$, and $\Gamma_2$ for the Hanoi Towers group $\HH$}
\label{f:hanoi}
\end{figure}
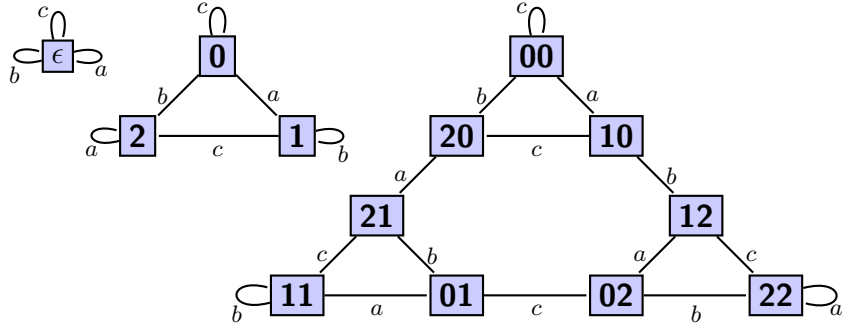

The group $\HH$ was introduced in~\cite{grigorchuk-s:hanoi-cr}. It models the well known Hanoi
Towers game on three pegs in such a way that the Schreier graph $\Gamma_n$ models the game for $n$
disks. It is the first example of a finitely generated branch group that admits a surjective
homomorphism onto the infinite dihedral group $D_\infty$ (note that branch groups can only have
virtually abelian proper quotients~\cite{grigorchuk:jibg}, and any finitely generated branch group
that admits a surjective homomorphism to an infinite virtually abelian group must map onto $\ZZ$ or
onto $D_\infty$~\cite{delzant-g:fa}).
\end{ex}

The \emph{boundary} $\boundary$ of the tree $X^*$ is the space of ends of the tree $X^*$. More
concretely, this is the space of all infinite rays
\[
 \boundary = \sset{x_1x_2x_3 \dots}{x_1,x_2,x_3,\dots \in X},
\]
that is, infinite paths without backtracking that start at the root. It has the structure of a
metric space (in fact, ultrametric space) with metric defined by $d(\xi,\zeta) = 1/2^{|\xi \wedge
\zeta|}$, where $\xi \wedge \zeta$ denotes the longest common prefix of the infinite rays $\xi$ and
$\zeta$, and $|\xi \wedge \zeta|$ denotes its length. Thus, the longer the common prefix the closer
the rays are. The induced topology is the product topology on $\prod_{i=1}^\infty X$, where the
finite space $X$ is given the discrete metric, implying that, topologically, the boundary
$\boundary$ is a Cantor set, and hence compact.

The action of any group of tree automorphisms $G \leq \Aut(X^*)$ naturally induces an action on the
boundary of the tree $X^*$. The action of any automorphism $g \in \Aut(X^*)$ on $\boundary$ is
given by~\eqref{e:gxw} with the understanding that $w$ in that formula now applies to rays in
$\boundary$, that is, to right-infinite words over $X$. If, for $n \geq 0$, we denote by
$\widetilde{\delta}_n: \boundary \to X^n$ the map that deletes the tail of any ray beyond the first
$n$-letters we obtain a sequence of $G$-equivariant maps. Thus we obtain the following diagram of
$G$-equivariant maps.
\begin{equation}\label{d:equivX}
\xymatrix{
 X^0 \ar@{<-}[r]^{\delta_0} & X^1 \ar@{<-}[r]^{\delta_1}  & X^2 \ar@{<-}[r]^{\delta_2} & \dots \\
 \boundary   \ar@{->}[u]^{\widetilde{\delta}_0} \ar@{->}[ur]^{\widetilde{\delta}_1} \ar@{->}[urr]^{\widetilde{\delta}_2} \ar@{..>}[urrr]
}
\end{equation}

Even if $G = \langle S \rangle$ acts \emph{level transitively} on the tree $X^*$ (transitively on
each level $X^n$ of the tree) and all Schreier graphs $\Gamma_n=\Gamma(G,S,X^n)$ are connected, the
Schreier graph $\Gamma(G,S,\boundary)$ of the action of $G$ on the tree boundary is not connected.
Indeed, since this graph is uncountable and the group $G$ is countable, each orbit of the action on
the boundary is countable and there must be uncountably many connected components (orbits) in the
graph $\Gamma(G,S,\boundary)$. Picking a connected component is equivalent to picking a point on
the boundary that represents it, that is, picking an infinite ray $\xi = x_1x_2x_3\dots \in
\boundary$. Choose such a ray $\xi$ and let $\Gamma= \Gamma_\xi = \Gamma(G,S,G\xi)$ be the Schreier
graph of the boundary action of $G$ on the orbit $G\xi = \{g(\xi) \mid g\in G\}$. We call the
Schreier graph $\Gamma=\Gamma_\xi$ the \emph{orbital Schreier graph} of $G$ at $\xi$. It is a
countable graph of degree $|S|$ and, since the restrictions of the maps $\widetilde{\delta}_n$, for
$n \geq 0$, to the orbit $G\xi$ are $G$-equivariant, the induced maps $\widetilde{\delta}_n: \Gamma
\to \Gamma_n$ are coverings. Therefore, we are precisely in the situation described by the
diagram~\eqref{d:equiv}. Moreover, we can now state a sufficient condition under which the spectra
of the sequence of finite graphs $\{\Gamma_n\}$ approximates the spectrum of $\Gamma$.

\begin{thm}[Bartholdi-Grigorchuk~\cite{bartholdi-g:spectrum}]\label{t:amenable-spectrum}
Let $G= \langle S \rangle \leq \Aut(X^*)$ be a finitely generated, self-similar, level-transitive
group of automorphisms of the rooted tree $X^*$ and let $\xi \in \boundary$ be a point on the tree
boundary. For $n \geq 0$, let $\Gamma_n = \Gamma(G,S,X^n)$ be the Schreier graph of the action of
$G$ on level $n$ of the tree and let $\Gamma= \Gamma_\xi = \Gamma(G,S,G\xi)$ be the orbital
Schreier graph of $G$ at $\xi$. If the action of $G$ on the orbit $G\xi$ is amenable, then
\[
 \overline{\bigcup_{n=0}^\infty \Sp(\Gamma_n)} = \Sp(\Gamma).
\]
\end{thm}

We recall the definition of an amenable action. The action of $G$ on $Y$ is \emph{amenable} if
there exists a normalized, finitely additive, $G$-invariant measure $\mu$ on all subsets of $Y$,
that is, there exists a function $\mu: 2^Y \to [0,1]$ such that
\begin{itemize}
\item (normalization) $\mu(Y)=1$,
\item (finite additivity) $\mu(A \sqcup B)= \mu(A)+\mu(B)$, for disjoint subsets $A,B \subseteq
    Y$,
\item ($G$-invariance) $\mu(A) = \mu(gA)$, for $g \in G$, $A \subseteq Y$.
\end{itemize}

For a finitely generated group $G=\langle S \rangle$ (with $S$ finite and symmetric, as usual)
acting transitively on on a set $Y$, the amenability of the action is equivalent to the amenability
of the Schreier graph $\Gamma=\Gamma(G,S,Y)$ of the action and one of the many equivalent ways to
define/characterize the amenability of $\Gamma$ is as follows. The graph $\Gamma$ is amenable if
and only if
\[
 \inf\sset{\frac{|\partial F|}{|F|}}{F \textup{ finite and nonempty set of vertices of } \Gamma} =0,
\]
where the \emph{boundary} $\partial F$ of the set $F$ is the set of vertices in $\Gamma$ that are
not in $F$ but have a neighbor in $F$, that is, $\partial F = \sset{ v \in \Gamma}
  {v \not\in F \textup{ and } sv \in F\textup{ for some } s \in S}$.

One sufficient condition for the amenability of the graph $\Gamma$ is obtained by looking at its
growth. Let $\Gamma$ be any connected graph of uniformly bounded degree. Choose any vertex $v_0 \in
\Gamma$ and, for $n \geq 0$, let $\gamma_{v_0}(n)$ be the number of vertices in $\Gamma$ at
combinatorial distance no greater than $n$ from $v_0$. If the growth of $\gamma_{v_0}(n)$ is
subexponential (that is, $\limsup_{n\to \infty} \sqrt[n]{\gamma_{v_0}(n)} =1$), then $\Gamma$ is an
amenable graph.

By definition, a group $G$ is amenable if its left regular action on itself is amenable. In such a
case, every action of $G$ is amenable and Theorem~\ref{t:amenable-spectrum} applies. The class of
amenable groups includes all finite and all solvable groups and is closed under taking subgroups,
homomorphic images, extensions, and directed unions. The smallest class of groups that contains all
finite and all abelian groups and is closed under taking subgroups, homomorphic images, extensions,
and directed unions is known as the class of \emph{elementary amenable} groups. There are amenable
groups that are not elementary amenable and many such examples came from the theory of self-similar
groups, starting with the first Grigorchuk group $\GG$. The amenability of this group was proved by
showing that it has subexponential (in fact intermediate, between polynomial and exponential)
growth~\cite{grigorchuk:gdegree}. Other examples of amenable but not elementary amenable groups
include Basilica group $\BB$~\cite{bartholdi-v:basilica}, Hanoi Towers group $\HH$, tangled
odometers group $\TT$, and many other automaton groups. See~\cite{bartholdi-k-n:bounded}
and~\cite{amir-a-v:linear} for useful sufficient conditions for amenability of automaton groups
based on random walk considerations and the notion of activity growth introduced by
Sidki~\cite{sidki:circuit}.
 
A large and interesting class of examples to which Theorem~\ref{t:amenable-spectrum} applies is the
class of contracting self-similar groups.

\begin{defn}
Let $G \leq \Aut(X^*)$ be a self-similar group of automorphisms of the rooted regular tree $X^*$.
The group $G$ is said to be \emph{contracting} if there exists a finite set $\mathcal{N}\subseteq
G$ such that, for every $g\in G$, there exists $n$ such that $g_v\in\mathcal{N}$, for all words
$v\in X^*$ of length at least $n$. The smallest set $\mathcal{N}$ satisfying this property is
called the \emph{nucleus} of the group.
\end{defn}

Since the growth of each orbital Schreier graph $\Gamma$ of a finitely generated, self-similar,
contracting group is polynomial~\cite{bartholdi-g:spectrum}, such a graph $\Gamma$ is amenable and,
therefore, its spectrum can be approximated by the spectra of the finite graphs in the sequence
$\{\Gamma_n\}$, as in Theorem~\ref{t:amenable-spectrum}. Note that it is not known yet whether all
finitely generated contracting groups are amenable.

\section{Iterated monodromy groups}\label{s:img}

The content of this section is not necessary in order to follow the rest of the survey, but it
provides excellent examples, motivation, and context for our considerations.

\subsection{Definition}
Let $\mathcal{M}$ be a path connected and locally path connected topological space, and let
$f:\mathcal{M}_1\to \MM$ be a finite degree covering map, where $\mathcal{M}_1$ is a subset of
$\MM$. The main examples for us are \emph{post-critically finite complex rational functions}.
Namely, a rational function $f\in\CC(z)$ is said to be post-critically finite if the forward orbit
$O_x=\{f^{\circ n}(x)\}_{n\ge 1}$ of every critical point $x$ of $f$ (seen as a self-map of the
Riemann sphere $\widehat\CC$) is finite. Let $P$ be the union of the forward orbits $O_x$, for all
critical points. Denote $\mathcal{M}=\widehat\CC \setminus P$ and
$\mathcal{M}_1=f^{-1}(\mathcal{M})$. Then $\mathcal{M}_1\subseteq \mathcal{M}$ and $f:\MM_1 \to
\MM$ is a finite degree covering map.

Let $t\in\mathcal{M}$, and consider the \emph{tree of preimages} $T_f$ whose set of vertices is the
disjoint union of the sets $f^{-n}(t)$, where $f^{-0}(t)=\{t\}$. We connect every vertex $v\in
f^{-n}(t)$ to the vertex $f(v)\in f^{-(n-1)}(t)$. We then obtain a tree rooted at $t$.

If $\gamma$ is a loop in $\mathcal{M}$ starting and ending at $t$ then, for every $v\in f^{-n}(t)$,
there exists a unique path $\gamma_v$ starting at $v$ such that $f^{\circ n}\circ\gamma_v=\gamma$.
Denote by $\gamma(v)$ the end of the path $\gamma_v$. Then $v\mapsto\gamma(v)$ is an automorphism
of the rooted tree $T_f$. We get in this way an action (called the \emph{iterated monodromy
action}) of the fundamental group $\pi_1(\mathcal{M}, t)$ on the rooted tree $T_f$. The quotient of
the fundamental group by the kernel of the action is called the \emph{iterated monodromy group} of
$f$, and is denoted $\img{f}$. In other words, $\img{f}$ is the group of all automorphisms of $T_f$
that are equal to a permutation of the form $v\mapsto\gamma(v)$ for some loop $\gamma$.

\subsection{Computation of $\img{f}$}

\begin{figure}
\centering
\includegraphics[width=95pt]{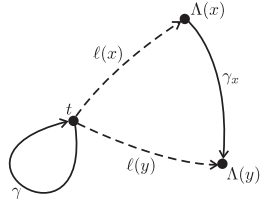}
\caption{Computation of $\img{f}$}
\label{f:recurn}
\end{figure}

Let $X$ be a finite alphabet of size $\deg f$, and let $\Lambda:X\to f^{-1}(t)$ be a bijection. For
every $x\in X$, choose a path $\ell(x)$ starting at $t$ and ending at $\Lambda(x)$. Let
$\gamma\in\pi_1(\MM, t)$. Denote by $\gamma_x$ the path starting at $\Lambda(x)$ such that
$f\circ\gamma_x=\gamma$, and let $\Lambda(y)$ be the end of $\gamma_x$. Then the paths $\ell(x)$,
$\gamma_x$, and $\ell(y)^{-1}$ form a loop, which we will denote  $\gamma|_x$ (see
Figure~\ref{f:recurn}).

\begin{prop}[Nekrashevych~\cite{nekrashevych:book-self-similar}]
Let $X$ be an alphabet in a bijection $\Lambda:X\longrightarrow f^{-1}(t)$. Let $\ell(x)$, $y$, and
$\gamma|_x$ be as above. Then $\Lambda$ can be extended to an isomorphism of rooted trees
$\Lambda:X^*\longrightarrow T_f$ that conjugates the iterated monodromy action of $\pi_1(\MM, t)$
on $T_f$ with the action on $X^*$ defined by the recursive rule:
\[\gamma(xv)=y\gamma|_x(v).\]
In particular, $\img{f}$ is a self-similar group.
\end{prop}

The self-similar action of $\img{f}$ on $X^*$ described in the last proposition is called the
\emph{standard action}. It depends on the choice of the connecting paths $\ell(x)$, for $x \in X$,
and the bijection $\Lambda:X\to f^{-1}(t)$. Changing the connecting paths amounts to
post-composition of the wreath recursion with an inner automorphism of the wreath product
$\Sym(X)\wr_X\img{f}$.

\begin{ex}[Basilica group $\BB=\img{z^2-1}$]
The polynomial $z^2-1$ is post-critically finite with $P=\{0, -1, \infty\}$. The fundamental group
of $\widehat\CC \setminus P$ is generated by two loops $a, b$ going around the punctures $0$ and
$-1$, respectively. With an appropriate choice of the connecting paths
(see~\cite[Subsection~5.2.2.]{nekrashevych:book-self-similar}), the wreath recursion for
$\img{z^2-1}$ is exactly the same as the one in Example~\ref{ex:basilica}. Thus, $\BB=\img{z^2-1}$.
\end{ex}

\begin{ex}[Tangled odometers group $\TT=\img{-\frac{z^3}{2}+\frac{3z}{2}}$]
The polynomial $f(z)=-z^3/2+3z/2$ has three critical points: $1$, $-1$, and $\infty$. All of them
are fixed points of $f$, hence $P=\{1, -1, \infty\}$, and the fundamental group of is generated by
loops around $1$ and $-1$. The corresponding iterated monodromy group is defined by the wreath
recursion~\eqref{e:tangled}, and this is the tangled odometers group $\TT$.
\end{ex}

\begin{ex}[Hanoi Towers group $\HH=\img{z^2-\frac{16}{27z}}$]
The iterated monodromy group of the rational function $z^2-16/(27z)$ is conjugate in $\Aut(X^*)$ to
the Hanoi Towers group $\HH$ (see~\cite{grigorchuk-s:standrews}).
\end{ex}

\begin{ex}[Dihedral group $D_\infty=\img{z^2-2}$ and binary odometer group $\ZZ=\img{z^2}$]
The iterated monodromy group of the polynomial $z^2-2$ is the dihedral group $D_\infty$ and of the
polynomial $z^2$ is the binary odometer group $\ZZ$ (infinite cyclic group) from
Example~\ref{ex:automata}.
\end{ex}

\subsection{Limit spaces of contracting self-similar groups}

Suppose that $G$ is a contracting self-similar group. Let $X^{-\omega}$ be the space of all
left-infinite sequences $\ldots x_2x_1$ of elements of $X$ with the direct product topology. We say
that two sequences $\ldots x_2x_1$ and $\ldots y_2y_1$ in $X^{-\omega}$ are \emph{asymptotically
equivalent} if there exists a sequence $\{g_k\}_{k=1}^\infty$ of elements in $G$, taking a finite
set of values, such that $g_k(x_k\ldots x_1)=y_k\ldots y_1$, for all $k\ge 1$. It is easy to see
that this is an equivalence relation. The \emph{limit space} of $G$ is the quotient of the
topological space $X^{-\omega}$ by the asymptotic equivalence relation. It is always a metrizable
space of finite topological dimension (if $G$ is contracting). Note that the asymptotic equivalence
relation is invariant with respect to the shift $\ldots x_2x_1\mapsto \ldots x_3x_2$. Consequently,
the shift induces a continuous self-map on the limit space of $G$. The obtained map is called the
\emph{limit dynamical system} of the group $G$.

\begin{thm}[Nekrashevych~\cite{nekrashevych:book-self-similar}]\label{t:correspondence}
Suppose that $f$ is a post-critically finite complex rational
function. Then $\img{f}$ is a contracting self-similar group with
respect to any standard action. The limit dynamical system of
$\img{f}$ is topologically conjugate to the restriction of $f$ onto
its Julia set.
\end{thm}

\begin{figure}
\begin{minipage}[c]{0.49\textwidth}
\begin{tabular}{c}
\includegraphics[width=100pt]{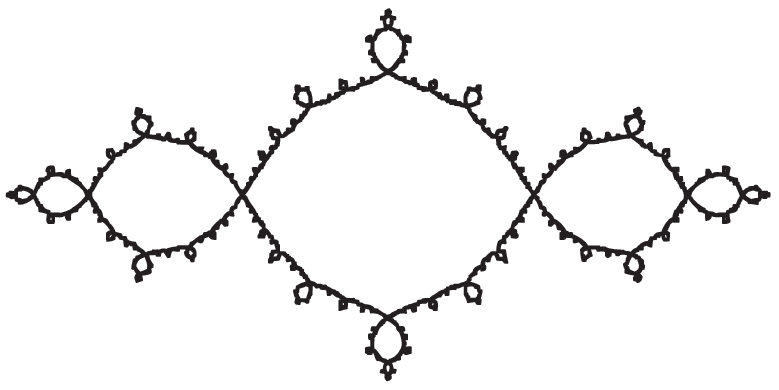}
\\
\includegraphics[width=125pt]{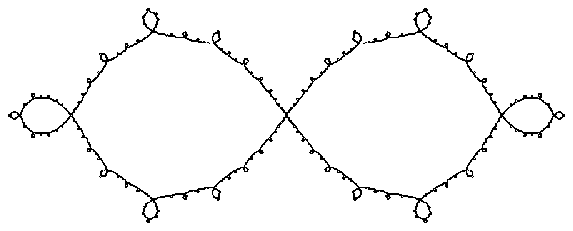}
\end{tabular}
\end{minipage}
\begin{minipage}[c]{0.49\textwidth}
\includegraphics[width=150pt]{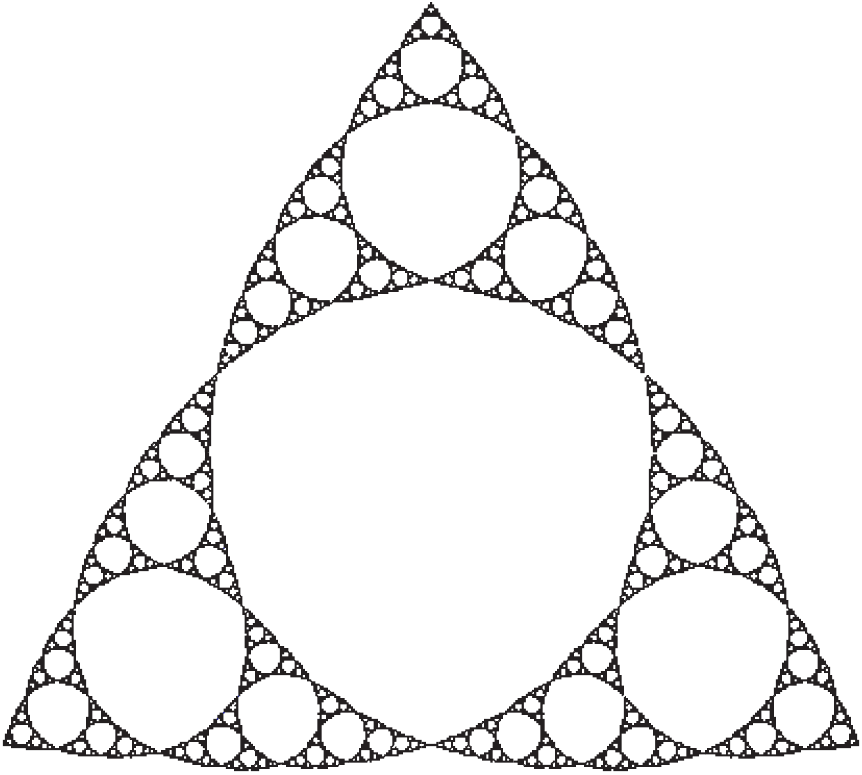}
\end{minipage}
\caption{Julia set of $z \mapsto z^2-1$ (top left), $z \mapsto -\frac{z^3}{2}+\frac{3z}{2}$ (bottom left), and $z \mapsto z^2-\frac{16}{27z}$ (right)}
\label{f:julia}
\end{figure}

The \emph{Julia set} of a complex rational function $f$ can be defined as the closure of the set of
points $c$ such that there exists $n$ such that $f^n(c)=c$ and $\left|\left.\frac d{dz}
f^n(z)\right|_{z=c}\right|>1$. The Julia sets of $z \mapsto z^2-1$, $z \mapsto
-\frac{z^3}{2}+\frac{3z}{2}$, and $z \mapsto z^2-\frac{16}{27z}$ are given in Figure~\ref{f:julia}.
Theorem~\ref{t:correspondence} provides context and explanation for the striking similarity between
the structure of the Schreier graphs of the Basilica group in Figure~\ref{f:basilica} and the
Basilica fractal in Figure~\ref{f:julia}, as well as between the structure of the Schreier graphs
of the Hanoi Towers group in Figure~\ref{f:hanoi} and the Sierpi\'nski gasket in
Figure~\ref{f:julia}.

\section{Relation to other operators and spectra}

\subsection{Hecke type operators}

Let $G = \langle S \rangle$, with $S$ finite and symmetric, be a finitely generated group and
$\lambda: G \to \unitary(\WW)$ a unitary representation of $G$ on a Hilbert space $\WW$. To each
element $\mathfrak{m} = \sum_{i=1}^n \alpha_i \cdot g_i$ of the group algebra $\CC[G]$ one can
associate the operator
\[
 \lambda(\mathfrak{m}) = \sum_{i=1}^n \alpha_i \lambda(g_i).
\]
In particular, we consider the \emph{Hecke type operator} $H_\lambda$ on the Hilbert space $\WW$
associated to the group algebra element $\mathfrak{h}=\displaystyle{\frac{1}{|S|}\sum_{s \in S} s}$
and given by
\[
  H_\lambda= \frac{1}{|S|}\sum_{s \in S} \lambda(s).
\]

\subsection{Koopman representation and Hecke type operators}

Let $G$ be a countable group acting on a measure space $(Y,\mu)$ by measure-preserving
transformations. The \emph{Koopman representation} $\pi$ is the unitary representation of $G$ on
the Hilbert space $L^2(Y, \mu)$ given by
\[ (\pi(g)f)(y) = f(g^{-1}y) \]
for $f \in L^2(Y, \mu)$ and $y \in Y$.

Let $G= \langle S \rangle \leq \Aut(X^*)$ be a finitely generated, self-similar, level-transitive
group of automorphisms of the rooted regular tree $X^*$. Note that the boundary $\boundary$, which
has the structure of a Cantor set $\prod_{i=1}^\infty X$, is a measure space with respect to the
product of uniform measures on $X$ (for the cylindrical set $uX^*$, we have $\mu(uX^*) =
\frac{1}{|X|^{|u|}}$). The group $G$ acts on $\boundary$ by measure-preserving transformations and
we may consider the Koopman representation $\pi$ of $G$ on $L^2(\boundary,\mu)$ and the associated
Hecke type operator $H_\pi$ on $L^2(\boundary,\mu)$, given by
\[
  H_\pi = \frac{1}{|S|}\sum_{s \in S} \pi(s).
\]
For every $n \geq 0$, we may also consider the representation $\pi_n$ on $L^2(X^n,\mu_n)$ on the
finite probability space $X^n$ with uniform probability measure $\mu_n$, corresponding to level $n$
of the tree, and the associated Hecke type operator
\[
 H_{\pi_n} = \frac{1}{|S|}\sum_{s \in S} \pi_n(s).
\]
Denote $\Sp(H_\pi) = \Sp(\pi)$ and $\Sp(H_{\pi_n}) = \Sp(\pi_n)$, for $n \geq 0$.

\begin{thm}[Bartholdi-Grigorchuk~\cite{bartholdi-g:spectrum}]
Let $G$ be a finitely generated, self-similar, level-transitive group of automorphisms of the
rooted regular tree $X^*$. Then
\[
 \Sp(\pi) = \overline{\bigcup_{n=0}^\infty \Sp(\pi_n)}.
\]
\end{thm}

Note that, unlike in Theorem~\ref{t:amenable-spectrum}, no additional requirements (such as
amenability of the action) are needed in the last result.

\subsection{Quasi-regular representations and Hecke type operators}

It is well known that every transitive left action of a group $G$ on any set $Y$ is equivalent to
the action of $G$ on the left coset space $G/P$, where $P=\Stab_G(y)$ is the stabilizer of the
point $y \in Y$ (since the action is transitive this point may be chosen arbitrarily). In fact,
Schreier graphs originate as the graphs of the action of groups on their coset spaces.

For a countable group $G$ and any subgroup $P \leq G$, the \emph{quasi-regular representation} is
the unitary representation $\rho_{G/P}$ of $G$ on the Hilbert space $\ell^2(G/P)$ given by
\[
 (\rho_{G/P}(g)f)(hP) = f(g^{-1}hP),
 \]
 for $f \in \ell^2(G/P)$  and $h \in G$. When $P$ is the trivial group we obtain the \emph{left regular representation} $\rho_G$ defined
 by
 \[
 (\rho_G(g)f)(h) = f(g^{-1}h),
\]
for $f \in \ell^2(G)$  and $h \in G$.

Let $G= \langle S \rangle \leq \Aut(X^*)$ be a finitely generated, self-similar, level-transitive
group of automorphisms of the rooted regular tree $X^*$ and let $\xi=x_1x_2x_3 \dots$ be a point on
the boundary $\boundary$. For $n \geq 0$, the point $x_1x_2\dots x_n$ is the unique point at level
$n$ on the ray $\xi$. Let
\begin{gather*}
 P_n = \Stab_G(x_1 x_2 \dots x_n), \text{ for } n \geq 0, \text{ and } \\
 P = \Stab_G(\xi).
\end{gather*}
Note that $\bigcap_{n=0}^\infty P_n = P_\xi$.

Denote by $\rho_n$ the quasi-regular representation $\rho_{G/P_n}$ corresponding to the subgroup
$P_n$ (thus, to the action of $G$ on level $n$ of the tree) and by $\rho_\xi$ the representation
$\rho_{G/P_\xi}$. We consider the Hecke type operator $H_{\rho_\xi}$ on $\ell^2(G/P_\xi)$
\[
  H_{\rho_\xi} = \frac{1}{|S|}\sum_{s \in S} \rho_\xi(s)
\]
and, for $n \geq 0$, the Hecke type operator
\[
 H_{\rho_n} = \frac{1}{|S|}\sum_{s \in S} \rho_n(s).
\]
Denote $\Sp(H_{\rho_\xi}) = \Sp(\rho_\xi)$ and $\Sp(H_{\rho_n}) = \Sp(\rho_n)$, for $n \geq 0$.

The following result extends Theorem~\ref{t:amenable-spectrum} and compares the Schreier spectrum
to the spectrum of the Hecke type operators $H_\pi$ and $H_{\rho_\xi}$ associated to the Koopman
representation $\pi$ and the quasi-regular representation $\rho_\xi$, respectively.

\begin{thm}[Bartholdi-Grigorchuk~\cite{bartholdi-g:spectrum}]\label{t:all-same}
(a) Let $G= \langle S \rangle \leq \Aut(X^*)$ be a finitely generated, self-similar,
level-transitive group of automorphisms of the rooted regular tree $X^*$ and let $\xi \in
\boundary$. Then, for $n \geq 0$,
\[
 \frac{1}{|S|} \Sp(\Gamma_n)  = \Sp(\rho_n) = \Sp(\pi_n)
\]
and
\[
 \frac{1}{|S|} \Sp(\Gamma_\xi) = \Sp(\rho_\xi) \subseteq \Sp(\pi).
\]

(b) If the action of $G$ on the orbit $G\xi$ is amenable, then
\[
 \frac{1}{|S|} \overline{\bigcup_{n=0}^\infty \Sp(\Gamma_n)} = \frac{1}{|S|} \Sp(\Gamma_\xi) = \Sp(\rho_\xi)  = \Sp(\pi).
\]

(c) If the group $P_\xi$ is amenable, then
\[
 \frac{1}{|S|} \Sp(\Gamma_\xi) = \Sp(\rho_\xi)  \subseteq \Sp(\rho_G),
\]
where $\rho_G$ is the left-regular representation of $G$ (and $\Sp(\rho_G)$ is the spectrum of the
corresponding Hecke type operator $H_{\rho_G}$).
\end{thm}

By part (b) in the last result, if the group $G$ is amenable, then all orbital Schreier graphs have
the same spectrum (there is no dependence on the choice of the point $\xi \in \boundary$, since the
representation $\pi$ does not depend on it). More generally, if all orbital Schreier graphs
$\Gamma_\xi$, for $\xi \in \boundary$ are amenable, as it happens in the case of contracting
self-similar groups, then they all have the same spectrum. Examples of nonamenable groups with
amenable orbital Schreier graphs $\Gamma_\xi$ were provided in~\cite{grigorchuk-n:nactions} (thus,
part (b) applies to some nonamenable groups).

We point out that part (b) is mistakenly stated in~\cite{bartholdi-g:spectrum} under the assumption
that either the action of $G$ on the orbit $G\xi$ is amenable or $P_\xi$ is amenable. The
assumption that $P_\xi$ is amenable only applies in part (c), and this part of
Theorem~\ref{t:all-same} follows from~~\cite[Proposition 3.5]{bartholdi-g:spectrum}.

\section{Method of computation}\label{s:method}

The method of computation of spectra, introduced in~\cite{bartholdi-g:spectrum} and further
implemented and refined
in~\cite{grigorchuk-z:l2,grigorchuk-s:hanoi-spectrum,grigorchuk-n:schur,grigorchuk-n-s:first-julia}
is based on the use of invariant sets of multidimensional rational maps and the Schur complement.
We will present the approach in the next two subsections, one addressing the global picture, and
the other the details.

\subsection{A global preview of the method}

Let $A$ be an operator for which we would like to calculate the spectrum. Include $A$ and the
entire pencil $\{ A(x) \mid x \in \CC \}$ with $A(x)=A-xI$ into a multidimensional pencil of operators
\[
 \{\ A^{(d)}(x_1,x_2,\dots,x_d) \mid x_1,\dots,x_d \in \CC \ \}
\]
such that
\[
 A(x) = A^{(d)}(x,x_2^{(0)},x_3^{(0)}, \dots, x_d^{(0)}),
\]
for some particular values $x_2^{(0)},x_3^{(0)},\dots,x_d^{(0)} \in \CC$. Define the joint spectrum
by
\[
 \Sp(A^{(d)}) =
  \sset{(x_1,x_2,\dots,x_d) \in \CC^d}{
 A^{(d)}(x_1,x_2,\dots,x_d) \textup{ is not invertible }}.
\]
Then
\[
 \Sp(A) = \Sp(A^{(d)}) \cap \ell,
\]
where $\ell$ is the line
\[
 \ell =
  \sset{(x_1,x_2,\dots,x_d) \in \CC^d}{x_2=x_2^{(0)}, x_3= x_3^{(0)},\dots,x_d = x_d^{(0)}}
\]
in the $d$-dimensional space $\CC^d$.

In the case of a self-adjoint operator $A$, which is always our case, we can use the field $\RR$
instead of $\CC$.

The problem naturally splits into three steps:
\begin{enumerate}
\item[(i)] Determine a suitable higher-dimensional pencil containing $\{ A(x) \mid x \in \RR \}$.

\item[(ii)] Determine the joint spectrum $\Sp(A^{(d)})$.

\item[(iii)] Determine the intersection $\Sp(A) = \Sp(A^{(d)}) \cap \ell$.
\end{enumerate}

In the examples that were successfully treated by this approach, the joint spectrum $\Sp(A^{(d)})$
is an invariant set under some rational $d$-dimensional map $F: \RR^d \to \RR^d$. Thus, in
practice, the step (ii) is understood as
\begin{enumerate}
\item[(ii)$'$] Determine the joint spectrum $\Sp(A^{(d)})$ as an $F$-invariant set for a suitable
    $d$-dimensional rational map $F: \RR^d \to \RR^d$.
\end{enumerate}

It may be somewhat counterintuitive why one should ``increase the dimension of the problem in order
to solve it'', but the method has worked well in situations were direct approaches have failed.
What happens is that the joint spectrum in $\RR^d$, corresponding to the $d$-fold pencil of
operators, is sometimes well behaved and easier to describe than the spectrum of the original
1-fold pencil. On the other hand, even when appropriate $A^{(d)}$ and $F$ are found, the structure
of the $F$-invariant set can be quite complicated and have the shape of a ``strange attractor''.

\subsection{More details}

Let $G = \langle S \rangle$ be an automaton group generated by the elements of the finite and
symmetric self-similar set $S$. For $n \geq 0$, the representations $\pi_n$ and $\rho_n$ are
equivalent and may be viewed as representations on the $|X|^n$-dimensional vector space
$\ell^2(X^n)$. The $|X|^n \times |X|^n$ adjacency matrix $A_n$ (the rows and the columns are
indexed by the words over $X$ of length $n$) of $\Gamma_n$ is given by
\[
 A_n = \sum_{s \in S} \pi_n(s).
\]
The $|X|^n \times |X|^n$ matrix $\pi_n(s)$ is given recursively, for $n > 0$, by blocks of size
$|X|^{n-1} \times |X|^{n-1}$
\begin{equation}\label{e:mwreath}
 \pi_n(s) = \left[ B_{y,x}(s)\right]_{y,x \in X}
\end{equation}
corresponding to the decomposition
\[ \ell^2(X^n) = \bigoplus_{x \in X} \ell^2(xX^{n-1}), \]
and the block $B_{y,x}(s)$ is given by
\[
 B_{yx}(s) =
 \begin{cases}
  \pi_{n-1}(s_x), & s(x)=y \\
  0         , & \text{otherwise}
 \end{cases}
\]
For $n=0$, the space $\ell^2(X^0)$ corresponding to the root of the tree is 1-dimensional and
$\pi_0(s)$ is the $1\times 1$ identity matrix $\pi_0(s)= [1]$. We call~\eqref{e:mwreath} the
\emph{matrix wreath recursion} of $S$ (it directly corresponds to the wreath recursion that defines
the generators $s \in S$).

From now on, we use the notation $s_n = \pi_n(s)$.

\begin{ex}
For the first Grigorchuk group $\GG$ the matrix wreath recursion gives
\[ a_0 = b_0 = c_0 = d_0 = [1] \]
and for $n >0$,
\[
 a_n = \begin{bmatrix} 0       & 1 \\ 1 & 0 \end{bmatrix} \qquad
 b_n = \begin{bmatrix} a_{n-1} & 0 \\ 0 & c_{n-1} \end{bmatrix} \qquad
 c_n = \begin{bmatrix} a_{n-1} & 0 \\ 0 & d_{n-1} \end{bmatrix} \qquad
 d_n = \begin{bmatrix} 1       & 0 \\ 0 & b_{n-1} \end{bmatrix},
\]
where, in each case, $0$ and $1$ denote the zero matrix and the identity matrix, respectively, of
appropriate size ($2^{n-1} \times 2^{n-1}$). Therefore, $A_0 = [4]$ and, for $n >0$,
\[
 A_n = \begin{bmatrix} 2a_{n-1}+1 & 1 \\ 1 & b_{n-1} + c_{n-1} + d_{n-1} \end{bmatrix}.
\]
\end{ex}

\begin{ex}
For the tangled odometers group $\TT$ the matrix wreath recursion gives
\[ a_0 = b_0 = a_0^{-1} = b_0^{-1} = [1] \]
and for $n \geq 0$,
\begin{alignat*}{2}
 a_{n+1} &= \begin{bmatrix} 0 & a_{n} & 0      \\ 1 & 0 & 0 \\ 0 & 0 & 1 \end{bmatrix} \qquad
 (a^{-1})_{n+1} &&=
   \begin{bmatrix} 0 & 1 & 0 \\ (a^{-1})_{n} & 0 & 0 \\ 0 & 0 & 1 \end{bmatrix} \\
 b_{n+1} &= \begin{bmatrix} 0 &  0 & b_{n-1}\\ 0 & 1 & 0 \\ 1 & 0 & 0 \end{bmatrix} \qquad
 (b^{-1})_{n+1} &&=
   \begin{bmatrix} 0 & 0 & 1 \\ 0 & 1 & 0 \\ (b^{-1})_n & 0 & 0 \end{bmatrix}.
\end{alignat*}
Therefore, $A_0 = [4]$ and, for $n \geq 0$,
\[
 A_{n+1} = \begin{bmatrix} 0 & 1+a_n & 1+b_n \\ 1+(a^{-1})_n & 2 & 0 \\ 1+(b^{-1})_n & 0 & 2\end{bmatrix}.
\]
\end{ex}

Once the recursive definition of the adjacency operator $A_n$ is established we consider the matrix
\[
 A_n(x) = A_n - xI = \left(\sum_{s \in S} s_n\right) - xI ,
\]
and more generally, a matrix of the form
\[
 A_n^{(d)}(x_1,\dots,x_d) = A_n - x_1I - \left(\sum_{i=2}^d x_i \cdot g_i \right) = \left(\sum_{s \in S} s_n\right) - x_1I - \left(\sum_{i=2}^d x_i \cdot g_i \right),
\]
for some auxiliary operators $g_2,\dots,g_d$. There is no known general approach how to choose
appropriate auxiliary operators. In practice, one needs to come up with good choices that make the
subsequent calculations feasible.

We then calculate, by using elementary column and row transformations and the Schur complement, the
determinant of $A_n^{(d)}$ in terms of the determinant of $A_{n-1}^{(d)}$ and obtain a recursive
expression of the form
\begin{equation}\label{e:det}
 \det(A_n^{(d)}(x_1,\dots,x_d)) = P_n(x_1,\dots,x_d) \det(A_{n-1}^{(d)}(F(x_1,\dots,x_d))),
\end{equation}
where $P_n(x_1,\dots,x_d)$ is a polynomial function and $F: \RR^d \to \RR^d$ is a rational function
in the variables $x_1,\dots,x_d$. Clearly, if the point $(x_1',\dots,x_d')$ is in the zero set of
$\det(A_{n-1}^{(d)}(x_1,\dots,x_d)$, then any point in $F^{-1}(x_1',\dots,x_d')$ is in the zero set
of $\det(A_n^{(d)}(x_1,\dots,x_d))$. Thus, describing the joint spectrum through iterations of the
recursion~\eqref{e:det} leads to iterations of the rational map $F$.

Understanding the structure of the zero sets of $\det(A_n^{(d)}(x_1,\dots,x_d))$, for $n \geq 0$,
and relating them to the zero sets of $\det(A_n(x))$ is accomplished, in the situations when we are
able to resolve this problem, by finding a function $\psi: \RR^d \to \RR$ and a polynomial function
$f: \RR \to \RR$ such that
\[
 \psi(F(x_1,\dots,x_d)) = f(\psi(x_1,\dots,x_d)),
\]
that is, by finding a semi-conjugacy from the $d$-dimensional rational function $F$ to a polynomial
function $f$ in a single variable. Since we have
\[
 \psi(F^{\circ m}(x_1,\dots,x_d)) = f^{\circ m}(\psi(x_1,\dots,x_d)),
\]
the iterations of $F$ are related to the iterations of $f$ and then the desired spectrum is
described through the iterations of the latter.

\section{Concrete examples and computation results}

In this section we present several concrete examples of calculations of spectra based on the method
suggested in the Section~\ref{s:method}. All groups in this section are amenable. By
Theorem~\ref{t:all-same}, the choice of the point on the boundary is irrelevant for the Schreier
spectrum and this is why no such choice is discussed in these examples.

One of the examples, the Hanoi Towers group $\HH$, leads to results on the Sierpi\'nski gasket. The
spectrum of Sierpi\'nski gasket goes back to the work of the physicists Rammal and
Toulouse~\cite{rammal-t:sierpinski}. It was turned into a mathematical framework by Fukushima and
Shima~\cite{fukushima-s:sierpinski}. Note that, in these works, the Sierpi\'nski gasket was
approximated by a sequence of graphs that are 4-regular (with the exception of the three corner
vertices, which have degree 2), while our approach yields an approximation through a different, but
related, sequence of 3-regular graphs. A method for spectra calculations in more general cases,
called \emph{spectral decimation}, was developed by Kumagai, Malozemov, Shima, Teplyaev, Strichartz
and
others~\cite{kumagai:decimation,malozemov:koch,shima:pcf,teplyaev:gasket,malozemov-t:self-similarity,strichartz:b-sierpinski}.
Connections with Julia sets are well-known, as for instance given by
Teplyaev~\cite{teplyaev:stochastic-survey}.

\subsection{The first Grigorchuk group $\GG$}
As was already mentioned, the method sketched above was introduced in~\cite{bartholdi-g:spectrum}
in order to compute the spectrum of the sequence of Schreier graphs $\{\Gamma_n\}$ and the boundary
Schreier graph $\Gamma$ for the case of the first Grigorchuk group $\GG$, as well as several other
examples, including the Gupta-Sidki 3-group~\cite{gupta-s:burnside}.

\begin{thm}[Bartholdi-Grigorchuk~\cite{bartholdi-g:spectrum}]\label{t:gg}
For $n \geq 1$, the spectrum of the graph $\Gamma_n$, as a set, has $2^n$ elements (thus, all
eigenvalues are distinct) and is equal to
\[
 \Sp(\Gamma_n) = \sset{1 \pm \sqrt{5+4\cos \frac{2k\pi}{2^n}}}{k=0,\dots,2^{n-1}} \setminus \{-2,0\}.
\]

The spectrum of $\Gamma$ (the Schreier spectrum of $\GG$), as a set, is equal to
\[
 \Sp(\Gamma) = [-2,0] \cup [2,4] .
\]
\end{thm}

\begin{rem}
There is a different way in which the spectrum of $\Gamma_n$ can be written. Namely, for $n \geq 2$,
\[
 \Sp(\Gamma_n) = \{4,2\} \cup \left ( 1 \pm \sqrt{5 \pm 2 \bigcup_{i=0}^{n-2} f^{-i}(0)} \right),
\]
where
\[
 f(x) = x^2 -2.
\]
Note that
\[
 f^{-k}(0) = \pm\sqrt{2 \pm \sqrt{2 \pm \sqrt{2\pm \dots \pm \sqrt{2}}} },
\]
where the root sign appears exactly $k$ times. The closure $\overline{\bigcup_{i=0}^\infty
f^{-i}(0)}$ is equal to the interval $[-2,2]$ and is the Julia set of the polynomial $f$,
Therefore,
\begin{align*}
 \Sp(\Gamma)
    &= \{4,2\} \cup \left(1 \pm \sqrt{5 \pm 2 \cdot [-2,2]}\right)
    = \{4,2\} \cup \left(1 \pm \sqrt{[1,9]}\right)  \\
    &= \{4,2\} \cup (1 \pm [1,3]) = [-2,0] \cup [2,4].
\end{align*}
\end{rem}

For the calculations in this example, we may use the 2-dimensional auxiliary pencil of operators
defined by
\[
 A_n^{(2)}(x,y) = a_n+b_n+c_n+d_n - (1+x)I + (y-1)a_n.
\]
The recursive formula for the determinant of $A_n(x,y)$ is, for $n \geq 2$,
\[
 \det(A_n^{(2)}(x,y)) = \left(x^2-4\right)^{2^{n-2}} \det(A_{n-1}^{(2)}(F(x,y))),
\]
where $F: \RR^2 \to \RR^2$ is given by
\[
 F(x,y) = \left(x-\frac{xy^2}{x^2-4},\frac{2y^2}{x^2-4}\right).
\]
The map $\psi: \RR^2 \to \RR$ that semi-conjugates $F$ to $f(x) = x^2 - 2$ is
\[
 \psi(x,y) = \frac{x^2-4 -y^2}{2y}.
\]
The $2$-dimensional joint spectrum of $A_n(x,y)$ is a family of hyperbolae and intersecting this
family with the line $y=1$ gives the desired spectrum.

The more general problem of determining the spectrum of the operator associated to any element of
the form $ta + ub + vc + wd$ in the group algebra $\RR[\GG]$ is considered
in~\cite{grigorchuk-l-n:all-coef}, where it is shown that, apart from few exceptions (such as the
case $u=v=w$ considered above), the spectrum is always a Cantor set.

\subsection{The Hanoi Towers group $\HH=\img{z^2 - \frac{16}{27z}}$ and Sierpi\'nski gasket}

\begin{thm}[Grigorchuk-{\v{S}}uni{\'c}~\cite{grigorchuk-s:hanoi-cr,grigorchuk-s:hanoi-spectrum}]\label{t:hh}
For $n \geq 1$, the spectrum of the graph $\Gamma_n$, as a set, has $3 \cdot 2^{n-1}-1$ elements and is equal to
\[
 \{3\} \ \cup \ \bigcup_{i=0}^{n-1} f^{-i}(0) \ \cup \
                \bigcup_{j=0}^{n-2} f^{-j}(-2),
\]
where
\[ f(x) = x^2 - x - 3. \]

The multiplicity of the $2^i$ eigenvalues in $f^{-i}(0)$, $i=0,\dots,n-1$ is
$a_{n-i}$, and the multiplicity of the $2^j$ eigenvalues in $f^{-j}(-2)$,
$j=0,\dots,n-2$, is $b_{n-j}$, where, for $m \geq 1$,
\[
 a_m = \frac{3^{m-1}+3}{2} \qqand b_m = \frac{3^{m-1}-1}{2}.
\]

The spectrum of $\Gamma$ (the Schreier spectrum of $\HH$), as a set, is equal to
\[
 \overline{\bigcup_{i=0}^\infty f^{-i}(0)} .
\]

It consists of a set of isolated points, the backward orbit $I= \bigcup_{i=0}^\infty f^{-i}(0)$ of
$0$ under $f$, and the set $J$ of accumulation points of $I$. The set $J$ is a Cantor set and is
the Julia set of the polynomial $f$.

The KNS spectral measure is concentrated on the union of the backward orbits
\[
 \left(\bigcup_{i=0}^\infty f^{-i}(0)\right) \cup \left(\bigcup_{i=0}^\infty f^{-i}(-2)\right).
\]
The KNS measure of each eigenvalue in $f^{-i}\{0,-2\}$, for $i=0,1,\dots$, is $\frac{1}{2 \cdot
3^{i+1}}$.

\end{thm}

\begin{rem}
The Kesten-von-Neumann-Serre measure (KNS measure for short) is the weak limit of the counting
spectral measures $\mu_n$ associated to the graph $\Gamma_n$, for $n \geq 0$ ($\mu_n(B) =
m_n(B)/|X|^n$, where $m_n(B)$ counts, including multiplicities, the eigenvalues of $\Gamma_n$ in
$B$.
\end{rem}

For the calculations in this example, the auxiliary pencil of operators used
in~\cite{grigorchuk-s:hanoi-spectrum} is 2-dimensional and given by
\[
 A_n^{(2)}(x,y) = a_n+b_n+c_n+ - xI + (y-1)d_n,
\]
where the block structure of $d_n$ is
\[
 d_n = \begin{bmatrix} 0 & 1 & 1 \\ 1 & 0 & 1 \\ 1 & 1& 0 \end{bmatrix}.
\]
The recursive formula for the determinant of $A_n^{(2)}(x,y)$ is, for $n \geq 2$,
\[
 \det(A_n^{(2)}(x,y)) = \left(x^2 - (1+y)^2\right)^{3^{n-2}} (x^2 - 1 + y - y^2)^{2 \cdot 3^{n-2}} \det(A_{n-1}^{(2)}(F(x,y))),
\]
where $F: \RR^2 \to \RR^2$ is given by
\[
 F(x,y) = \left(x+\frac{2y^2(-x^2+x+y^2)}{(x-1-y)(x^2-1+y-y^2)},\frac{y^2(x-1+y)}{(x-1-y)(x^2-1+y-y^2)}\right).
\]
The map $\psi: \RR^2 \to \RR$ that semi-conjugates $F$ to $f(x) = x^2 - x - 3$ is
\[
 \psi(x,y) = \frac{x^2-1-xy-2y^2}{y}.
\]

\subsection{The Tangled Odometers Group $\TT = \img{-\frac{z^3}{2}+\frac{3z}{2}}$ and the first Julia set}

\begin{thm}[Grigorchuk-Nekrashevych-{\v{S}}uni{\'c}~\cite{grigorchuk-n-s:first-julia}]\label{t:tt}
For $n \geq 0$, the spectrum of the graph $\Gamma_n$, as a set, has $2^{n+1}-1$
elements and is equal to
\[
 \{4\} \ \cup \ \bigcup_{i=0}^{n-1} f^{-i}(2) \ \cup \
                \bigcup_{j=0}^{n-1} f^{-j}(-2),
\]
where
\[ f(x) = x^2 - 2x - 4. \]

The multiplicity of the $2^i$ eigenvalues in $f^{-i}(2)$, $i=0,\dots,n-1$ is
$3^{n-1-i}$, the multiplicity of the $2^j$ eigenvalues in $f^{-j}(-2)$,
$j=0,\dots,n-1$, is 1, and the multiplicity of the eigenvalue 4 is 1.

The spectrum of $\Gamma$ (the Schreier spectrum of $\TT$), as a set, is equal to
\[
 \overline{\bigcup_{i=0}^\infty f^{-i}(2)}.
\]

It consists of a set of isolated points, the backward orbit $I= \bigcup_{i=0}^\infty f^{-i}(2)$ of
$2$ under $f$, and the set $J$ of accumulation points of $I$. The set $J$ is a Cantor set and is
the Julia set of the polynomial $f$.

The KNS spectral measure is concentrated on the backward orbit
\[
 I = \bigcup_{i=0}^\infty f^{-i}(2)
\]
of $f$. The KNS measure of each eigenvalue in $f^{-i}\{2\}$, for $i=0,1,\dots$, is
$\frac{1}{3^{i+1}}$.
\end{thm}

For the calculations in this example, the auxiliary pencil of operators used
in~\cite{grigorchuk-n-s:first-julia} is 3-dimensional and given by
\[
 A_n^{(3)}(x,y,z) = a_n+b_n+a_n^{-1}+b_n^{-1} - xc_n - (z+2)d_n + (y-1)g_n,
\]
where the block structure of $c_n$, $d_n$ and $g_n$ is
\[
 c_{n} =
   \begin{bmatrix}
     1 & 0 & 0 \\
     0 & 0 & 0 \\
     0 & 0 & 0
   \end{bmatrix}, \qquad
 d_{n+1} =
   \begin{bmatrix}
     0 & 0 & 0 \\
     0 & 1 & 0 \\
     0 & 0 & 1
   \end{bmatrix}, \qquad
 e_{n+1} =
   \begin{bmatrix}
     0 & 1 & 1 \\
     1 & 0 & 0 \\
     1 & 0 & 0
   \end{bmatrix}.
\]

\subsection{Lamplighter group $\LL_2 = \ZZ \ltimes \oplus_\ZZ \ZZ/2\ZZ$}

\begin{thm}[Grigorchuk-{\.Z}uk~\cite{grigorchuk-z:l2}]\label{t:ll}
For $n \geq 0$, the spectrum of the graph $\Gamma_n$, as a set, is equal to
\[
 \Sp(\Gamma_n) = \{4\} \cup
  \sset{4\cos\frac{p}{q}\pi}{1\leq p<q \leq n+1 \textup{ and } p \textup{ and } q \textup{relatively prime}} .
\]

The multiplicity of the eigenvalue $4\cos\frac{p}{q}\pi$, for $1 \leq p < q \leq n+1$, and $p$ and
$q$ relatively prime is equal to
\[
   \frac{2^n-2^{\mmod(n,q)}}{2^q-1} + 1_{\{q \textup{ divides }n+1\}},
\]
where $\mmod(n,q)$ is the remainder obtained when $n$ is divided by $q$, and $1_{\{q \textup{
divides }n+1\}}$ is the indicator function equal to 1 when $q$ divides $n+1$ and to 0 otherwise.
The multiplicity of the eigenvalue 4 is 1.

The spectrum of $\Gamma$ (the Schreier spectrum of $\LL_2$), as a set, is equal to
\[
 \Sp(\Gamma) = [-4,4].
\]

The KNS spectral measure is discrete and, for the eigenvalue $4 \cos \frac{p}{q} \pi$, with $1 \leq
p < q$ and $p$ and $q$ relatively prime, is equal to $\frac{1}{2^q-1}$.
\end{thm}

The above result has several interesting corollaries. First, note that there exist an infinite ray
$\zeta \in \partial X^*$ for which the corresponding parabolic subgroup $P_\zeta =
\St_{\LL_2}(\zeta)$ is trivial~\cite{grigorchuk-z:l2} (in fact, this is true for all infinite rays
that are not eventually periodic~\cite{nekrashevych-p:scale,grigorchuk-k:l2subgroups}). For such a
ray $\zeta$, the Schreier graph $\Gamma_\zeta=\Gamma(\LL_2,P_\zeta,S)$ and the Cayley graph
$\Gamma(\LL_2,S)$ are isomorphic. The calculation of the spectrum of $\LL_2$ led to a
counterexample of the Strong Atiyah Conjecture. The Strong Atiyah Conjecture states that if $\MM$
is a closed Riemannian manifold with fundamental group $G$, then its $L^2$-Betti numbers come from
the following subgroup of the additive group of rational numbers
\[
 \fin^{-1}(G) = \left\langle \sset{\frac{1}{|H|}}{H \textup{ a finite subgroup of }G} \right\rangle \leq \QQ.
\]
This is contradicted by the following result.

\begin{thm}[Grigorchuk, Linnell, Schick, \.{Z}uk \cite{grigorchuk-lsz:atiyah}]
There exists a closed Riemannion 7-dimensional manifold $\MM$ such that all finite groups in its
fundamental group $G$ are elementary 2-abelian, $\fin^{-1}(G)=\ZZ[\frac{1}{2}]$, but its third
$L^2$-Betti number is $\beta_{(2)}^{3}(\MM) = \frac{1}{3}$.
\end{thm}

Note that other versions of Atiyah Conjecture were later also disproved by using examples based on
lamplighter-like groups~\cite{austin:atiyah,lehner-w:atiyah}.

\subsection{Basilica group $\BB = \img{z^2-1}$ and $\img{z^2+i}$}

We do not have complete results for these two examples, but some progress was achieved.

The Schreier spectrum of Basilica group $\BB$ was considered in~\cite{grigorchuk-z:basilica2},
using the auxiliary 2-dimensional pencil of operators given by
\[
 A_n^{(2)}(x,y) = a_n+a_n^{-1}+ y(b_n^{-1}+b_n^{-1}) - xI.
\]
Partial results were also obtained by Rogers and Teplyaev by using the spectral decimation
method~\cite{rogers-t:basilica}.

The group $K =\img{z^2+i}$ of binary tree automorphisms is generated by three involutions defined
by the wreath recursion
\[
 a = (01)(e,e) \qquad  b = ()(a,c) \qquad  c = ()(b,e).
\]
The Schreier spectrum of $\img{z^2+i}$ was considered in~\cite{grigorchuk-s-s:z2i}, using the
auxiliary 3-dimensional pencil of operators given by
\[
 A_n^{(3)}(x,y,z) = a_n + yb_n + zc_n - xI.
\]
In both cases, the corresponding multi-dimensional map $F: \RR^d \to \RR^2$ was found, but the
shape of the corresponding $F$-invariant subset (that is, the joint spectrum) is unknown.

\section{Laplacians on the limit fractals}\label{s:laplacian}

For some contracting self-similar groups $G$,
the Hecke type operators $H_{\pi_n}$, when appropriately rescaled, converge
to a well defined \emph{Laplacian} on the limit space. The process of
finding the rescaling coefficient and proving existence of the limit
Laplacian has much in common with the process of computing the spectra
of operators $H_{\pi_n}$, as described in Section~\ref{s:method}. A
general theory, working for all contracting groups is still missing,
but many interesting examples can be analyzed.

The technique in the known examples is based on the theory of Dirichlet forms on self-similar sets,
see~\cite{kigami:book}. A connection of this theory with self-similar groups, and the examples
described in this section are discussed in more detail in~\cite{nektepl:analysis}.

Let $G$ be a self-similar group generated by a finite symmetric set $S$ and, for $n \geq 0$, let
$L_n=1-H_{\pi_n}$ be the corresponding Laplacian on the Schreier graph $\Gamma_n$. Let
$\mathcal{E}_n$ be the quadratic form with matrix $L_n$, that is, the form given by
$\mathcal{E}_n(x, y)=\langle L_nx, y\rangle$.

Choose a letter $x_0\in X$, and consider for every $n\ge 1$ the subset $V_n=x_0^{-\omega}X^n$ of
the space $X^{-\omega}$ encoding the limit space of $G$. We have $V_n\subseteq V_{n+1}$, and we
naturally identify $V_n$ with $X^n$ by the bijection $v\mapsto x_0^{-\omega}v$. We also consider
$\mathcal{E}_n$ as a form on $\ell^2(V_n)=\ell^2(X^n)$.

The \emph{trace} $\mathcal{E}_{n+1}'$ of $\mathcal{E}_{n+1}$ on $V_n$ is the
quadratic form $\mathcal{E}$ such that for $f\in\ell^2(V_n)$ the
value of $\mathcal{E}(f, f)$ is equal to the infimum of values of
$\mathcal{E}_{n+1}(g, g)$ over all functions $g\in\ell^2(V_{n+1})$
such that $g|_{V_n}=f$.

The matrix of $\mathcal{E}_{n+1}'$ is found as the \emph{Schur
  complement} of the matrix $L_{n+1}$ of $\mathcal{E}_{n+1}$. Namely,
decompose the matrix $L_{n+1}$ into the block form $\left[\begin{array}{cc} A & B\\ C &
D\end{array}\right]$ according to the decomposition of $\ell^2(X^{n+1})$ into the direct sum
$\ell^2(x_0X^n)\oplus\ell^2((X\setminus\{x_0\})X^n)$ (so that $A$, $B$, $C$, and $D$ are of sizes
$k^n\times k^n$, $k^n\times (k-1)k^n$, $(k-1)k^n\times k^n$, and $(k-1)k^n\times (k-1)k^n$,
respectively, where $k=|X|$ is the size of the alphabet). Then the matrix of $\mathcal{E}_{n+1}'$
is $A-BD^{-1}C$.

Let us consider some examples. Let $G=\img{z^2-1}$ be the Basilica group. Consider the Laplacian
$1-\alpha(a+a^{-1})-\beta(b+b^{-1})$, and the corresponding Dirichlet forms $\mathcal{E}_n$. Then
it follows from the recursive definition of the generators $a$ and $b$ that the decomposition of
$L_{n+1}$ into blocks $\left[\begin{array}{cc} A & B\\ C&
    D\end{array}\right]$ (for $x_0=1$) is
\[\left[\begin{array}{cc} 1-\beta(a+a^{-1}) & -\alpha(1+b^{-1})\\
-\alpha (1+b) & 1-2\beta\end{array}\right],\]
hence the matrix of $\mathcal{E}_{n+1, x_0}$ is
$\left(1-\frac\alpha 2\right)-\left(\beta(a+a^{-1})+\frac\alpha
  2\right)$. Consequently, if we take $\alpha=\frac{2-\sqrt{2}}2$, and
$\beta=\frac{\sqrt{2}-1}2$, then we have $\mathcal{E}_{n+1}'=\lambda\mathcal{E}_n$ for
$\lambda=\frac{1}{\sqrt{2}}$. It follows then from the general theory,
see~\cite{kigami:resistanceforms}, that the forms $\lambda^{-n}\mathcal{E}_n$ converge to a
Laplacian on the limit space of $G$, that is, on the Julia set of $z^2-1$.

In some cases one needs to take slightly bigger sets $V_n$. For
example, consider the Hanoi Towers group $\mathcal{H}$. Let $V_n$ be
the set of sequences of the form $0^{-\infty}X^n$, $1^{-\infty}X^n$, and
$2^{-\infty}X^n$. Let $a=(01)(e, e, a)$, $b=(02)(e, b, e)$, and
$c=(12)(c, e, e)$, and consider, for positive real numbers $x$, $y$,
the form $\mathcal{E}_n$ on $\ell^2(V_n)$ given by the matrix
\[\left[\begin{array}{ccc} y(1-a)-2x & -x & -x\\ -x & y(1-b)-2x & -x\\
    -x & -x & y(1-c)-2x\end{array}\right]\]
with respect to the decomposition
$\ell^2(V_n)=\ell^2(0^{-\omega}X^n)\oplus\ell^2(1^{-\omega}X^n)\oplus\ell^2(2^{-\omega}X^n)$,
  where $a, b, c$ act on the corresponding subspaces
  $\ell^2(x^{-\omega}X^n)$ using the representation $\pi_n$ (after we
  identify $x^{-\omega}X^n$ with $X^n$ in the natural way).

Then a direct computation using the recursive definition of the generators $a, b, c$, and the Schur
complement shows that trace of $\mathcal{E}_{n+1}$ on $V_n$ is given by the same matrix where $(x,
y)$ is replaced by $\left(\frac{3}{5+3x/y}x, y\right)$. Passing to the limit $y\to\infty$, and
restricting to functions on which the limit of the quadratic form is finite (which will correspond
to identifying sequences $\ldots x^{-\omega}v$ representing the same points of the limit space), we
get rescaling $x\mapsto \frac 35x$, hence convergence of $(5/3)^n\mathcal{E}_n$ to a Laplacian on
the limit space of $\mathcal{H}$, which is the Sierpi\'nski gasket.

 ------------------------------------------------------------------------

\subsection*{Acknowledgment}
Many thanks to Christoph Bandt for his generous help and valuable input.

\newcommand{\etalchar}[1]{$^{#1}$}
\def\cprime{$'$}


\end{document}